\documentclass[11pt,oneside]{amsart}
\pdfoutput=1

\usepackage{tikz-cd}

\usepackage{tcolorbox}

\usepackage[arrow,curve,matrix,tips,2cell]{xy}

\usepackage[hmargin=2cm,vmargin=1.5cm]{geometry}

\usepackage{bm}

\usepackage{fancyhdr}

\usepackage{amsmath}
\usepackage{amscd}
\usepackage{amssymb}
\usepackage{latexsym}

\usepackage{graphicx}

\usepackage{hyperref}
\usepackage{cleveref}
\usepackage{comment}

\newtheorem{theorem}{Theorem}[section]
\newtheorem*{theorem*}{Theorem}
\newtheorem{lemma}[theorem]{Lemma}
\newtheorem*{lemma*}{Lemma}
\newtheorem{corollary}[theorem]{Corollary}
\newtheorem{proposition}[theorem]{Proposition}

\newtheorem{remark}[theorem]{Remark}
\newtheorem{definition}[theorem]{Definition}
\newtheorem*{definition*}{Definition}

\newtheorem{question}[theorem]{Question}
\newtheorem*{question*}{Question}

\newtheorem{example}[theorem]{Example}
\newtheorem{examples}[theorem]{Examples}

\makeatletter
\def\revddots{\mathinner{\mkern1mu\raise\p@
\vbox{\kern7\p@\hbox{.}}\mkern2mu
\raise4\p@\hbox{.}\mkern2mu\raise7\p@\hbox{.}\mkern1mu}}
\makeatother 
\newcommand{\bgl}{\begin{equation}} 
\newcommand{\egl}{\end{equation}}
\newcommand{\bgloz}{\begin{equation*}} 
\newcommand{\egloz}{\end{equation*}}
\newcommand{\bgln}{\begin{eqnarray}} 
\newcommand{\egln}{\end{eqnarray}}
\newcommand{\bglnoz}{\begin{eqnarray*}} 
\newcommand{\eglnoz}{\end{eqnarray*}}
\newcommand{\btheo}{\begin{theorem}}
\newcommand{\etheo}{\end{theorem}}
\newcommand{\btheooz}{\begin{theorem*}}
\newcommand{\etheooz}{\end{theorem*}}

\newcommand{\blemma}{\begin{lemma}}
\newcommand{\elemma}{\end{lemma}}
\newcommand{\blemmaoz}{\begin{lemma*}}
\newcommand{\elemmaoz}{\end{lemma*}}
\newcommand{\bproof}{\begin{proof}}
\newcommand{\eproof}{\end{proof}}
\newcommand{\bbew}{\begin{beweis}}
\newcommand{\ebew}{\end{beweis}}
\newcommand{\bremark}{\begin{remark}\em}
\newcommand{\eremark}{\end{remark}}
\newcommand{\bdefin}{\begin{definition}}
\newcommand{\edefin}{\end{definition}}
\newcommand{\bdefinoz}{\begin{definition*}}
\newcommand{\edefinoz}{\end{definition*}}
\newcommand{\bex}{\begin{example}}
\newcommand{\eex}{\end{example}}
\newcommand{\bexs}{\begin{examples}}
\newcommand{\eexs}{\end{examples}}

\newcommand{\bprop}{\begin{proposition}}
\newcommand{\eprop}{\end{proposition}}
\newcommand{\bcor}{\begin{corollary}}
\newcommand{\ecor}{\end{corollary}}
\newcommand{\bfa}{\begin{cases}} 
\newcommand{\efa}{\end{cases}}

\newcommand{\bquestion}{\begin{question}}
\newcommand{\equestion}{\end{question}}
\newcommand{\bquestionoz}{\begin{question*}}
\newcommand{\equestionoz}{\end{question*}}
\newcommand{\n}{\par\noindent}

\newcommand{\cF}{\mathcal F}
\newcommand{\cG}{\mathcal G}

\newcommand{\cI}{\mathcal I}
\newcommand{\cJ}{\mathcal J}

\newcommand{\cM}{\mathcal M}

\newcommand{\cO}{\mathcal O}
\newcommand{\cP}{\mathcal P}
\newcommand{\cQ}{\mathcal Q}

\newcommand{\cS}{\mathcal S}

\def\Az{\mathbb{A}}
\def\Cz{\mathbb{C}}

\def\Kz{\mathbb{K}}

\def\Nz{\mathbb{N}}

\def\Qz{\mathbb{Q}}
\def\Rz{\mathbb{R}}
\def\Tz{\mathbb{T}}
\def\Zz{\mathbb{Z}}

\def\1z{\mathbb{1}}
\newcommand{\fC}{\mathfrak C}

\newcommand{\fP}{\mathfrak P}

\newcommand{\fX}{\mathfrak X}

\newcommand{\an}[1]{``#1''} 
\newcommand{\ti}{\tilde}

\newcommand{\lori}{\longrightarrow}

\newcommand{\ma}{\mapsto} 
\newcommand{\onto}{\twoheadrightarrow} 
\newcommand{\into}{\hookrightarrow} 
\newcommand{\Rarr}{\Rightarrow} 

\newcommand{\ve}{\varepsilon}

\def\SEMI{\mbox{$\times\kern-2pt\vrule height5pt width.6pt \kern3pt $}}

\newcommand{\Spec}{{\rm Spec\,}} 
\newcommand{\ind}{{\rm ind\,}}

\newcommand{\id}{{\rm id}}

\newcommand{\rk}{{\rm rk}}
\newcommand{\img}{{\rm im\,}}
\renewcommand{\ker}{{\rm ker}\,}

\newcommand{\reg}{^\times} 
\newcommand{\abs}[1]{\lvert#1\rvert} 
\newcommand{\defeq}{\mathrel{:=}} 

\newcommand{\dop}{\text{: }} 
\newcommand{\falls}{\text{ if }} 
\newcommand{\sonst}{\text{ else}} 
\newcommand{\ilim}{\varinjlim} 
\newcommand{\dotcup}{\ensuremath{\mathaccent\cdot\cup}} 
\newcommand{\supp}{{\rm supp}} 
\newcommand{\lge}{\left\{} 
\newcommand{\rge}{\right\}} 
\newcommand{\lru}{\left(} 
\newcommand{\rru}{\right)} 
\newcommand{\lsp}{\left\langle} 
\newcommand{\rsp}{\right\rangle} 
\newcommand{\rukl}[1]{\lru #1 \rru} 
\newcommand{\gekl}[1]{\lge #1 \rge} 
\newcommand{\spkl}[1]{\lsp #1 \rsp} 
\newcommand{\menge}[2]{\gekl{ #1 \dop #2 }} 
\newcommand{\mfa}{\mathfrak{a}}
\newcommand{\mfb}{\mathfrak{b}}
\newcommand{\mfm}{\mathfrak{m}}
\newcommand{\mfn}{\mathfrak{n}}
\newcommand{\mfp}{\mathfrak{p}}
\renewcommand{\a}{\mathfrak{a}}

\renewcommand{\k}{\mathfrak{k}}

\newcommand{\p}{\mathfrak{p}}
\newcommand{\q}{\mathfrak{q}}
\newcommand{\m}{\mathfrak{m}}
\renewcommand{\n}{\mathfrak{n}}

\newcommand{\Cl}{\textup{Cl}}
\newcommand{\Gal}{\textup{Gal}}
\newcommand{\Prim}{\textup{Prim}}

\newcommand{\sign}{\textup{sign}}

\newcommand{\acts}{\curvearrowright}

\def\clmg{{\rm Cl}_\m^{\bar{\Gamma}}}
\def\clnl{{\rm Cl}_\n^{\bar{\Lambda}}}
\def\rmg{R_{\m,\Gamma}}

\def\tor{{\rm tor}}
\def\Kmg{K(\m)^{\bar{\Gamma}}}
\def\Km{K(\m)}
\def\Knl{K(\n)^{\bar{\Lambda}}}

\newcommand{\ord}{\textup{ord}}

\newcommand{\im}{\textup{im}}

\begin{document}

\title{On K-theoretic invariants of semigroup C*-algebras from actions of congruence monoids}

\thispagestyle{fancy}

\author{Chris Bruce}
\address[Chris Bruce]{School of Mathematical Sciences, Queen Mary University of London, Mile End Road, E1 4NS London, United Kingdom, and
	School of Mathematics and Statistics, University of Glasgow, University Place, Glasgow G12 8QQ, United Kingdom}
\email{Chris.Bruce@glasgow.ac.uk}

\author{Xin Li}
\address[Xin Li]{School of Mathematics and Statistics, University of Glasgow, University Place, Glasgow G12 8QQ, United Kingdom}
\email{xin.li@glasgow.ac.uk}

\thanks{The first-named author was supported by the Natural Sciences and Engineering Research Council of Canada through an Alexander Graham Bell CGS-D award.}

\begin{abstract}
We study semigroup C*-algebras of semigroups associated with number fields and initial data arising naturally from class field theory. These semigroup C*-algebras turn out to have an interesting C*-algebraic structure, giving access to many new examples of classifiable C*-algebras and exhibiting phenomena which have not appeared before. Moreover, using K-theoretic invariants, we investigate how much information about the initial number-theoretic data is encoded in our semigroup C*-algebras. 
\end{abstract}

\maketitle

\setlength{\parindent}{0cm} \setlength{\parskip}{0.5cm}

\section{Introduction}

Semigroup C*-algebras are C*-algebras generated by left regular representations of left-cancellative semigroups. They form a natural example class of C*-algebras and have been studied in various contexts for several families of semigroups, for instance positive cones in totally ordered groups \cite{Co,Dou,Mur}, semigroups naturally arising in combinatorial or geometric group theory \cite{CL02,CL07,Sp12}, semigroups given by presentations \cite{LOS}, or semigroups of number-theoretic origin \cite{LacaRae,CDL,EL,Li14,Li16_2}. It was this last class of examples which has triggered many of the recent developments in semigroup C*-algebras (see \cite{CELY} and the references therein for an overview).

While semigroup C*-algebras for the full $ax+b$-semigroups over rings of algebraic integers have been studied in \cite{CDL,EL,Li14,Li16_2}, the first-named author considered a much more general class of semigroups in \cite{Bru,Bru2} and showed that, while providing a rich source of new examples, they allow for a similar analysis as in the full $ax+b$ case. This generalization, which is very natural from the number-theoretic point of view, proceeds as follows: Given a number field $K$ with ring of algebraic integers $R$, a modulus $\m$ for $K$, and a group $\Gamma$ of residues modulo $\m$, define the associated congruence monoid $R_{\m,\Gamma}$ as the multiplicative submonoid of elements in $R$ that are relatively prime to $\m$ and reduce to an element of $\Gamma$ modulo $\m$, and form the semi-direct product $R \rtimes R_{\m,\Gamma}$ where $R_{\m,\Gamma}$ acts on $R$ by multiplication. In this way, we obtain a generalization of the construction of $ax+b$-semigroups (the latter corresponding to the case of trivial $\m$ and $\Gamma$). The data $(\m,\Gamma)$ canonically gives rise to a class field $K(\m)^{\bar{\Gamma}}$ of $K$, which, in the case of trivial $\m$ and $\Gamma$, is simply the Hilbert class field of $K$. The semigroup C*-algebras attached to these generalized $ax+b$-semigroups turn out to be very interesting from the C*-algebraic perspective. In particular, their K-theoretic invariants are very rich and exhibit interesting new phenomena which have not appeared before.

Our goal in this paper is twofold. The first goal is a careful analysis of the semigroup C*-algebras $C^*_{\lambda}(R \rtimes R_{\m,\Gamma})$ and their K-theoretic invariants. Building on this, the second goal is to address the following natural question:
\setlength{\parindent}{0cm} \setlength{\parskip}{0cm}
\begin{center}
How much of the initial number field $K$ and the class field $K(\m)^{\bar{\Gamma}}$
\setlength{\parindent}{0cm} \setlength{\parskip}{0cm}

does our semigroup C*-algebra $C^*_{\lambda}(R \rtimes R_{\m,\Gamma})$ remember? 
\end{center}
\setlength{\parindent}{0cm} \setlength{\parskip}{0.5cm}

More precisely, our goal is to extract information about $K$ and $K(\m)^{\bar{\Gamma}}$ from K-theoretic invariants of $C^*_{\lambda}(R \rtimes R_{\m,\Gamma})$. The idea of using K-theory to extract information goes back to the classification programme of C*-algebras and has already proven to be fruitful in the case of $ax+b$-semigroups  \cite{Li14,Li16_2}. Our main result reads as follows:
\btheooz[see Theorem~\ref{thm:reconstruction}]
Suppose that $K$ and $L$ are number fields with rings of algebraic integers $R$ and $S$. Let $\m$ and $\n$ be moduli for $K$ and $L$, and let $\Gamma$ and $\Lambda$ be subgroups of $(R/\m)^*$ and $(S/\n)^*$, respectively. Suppose that there is an isomorphism $C_\lambda^*(R\rtimes R_{\m,\Gamma})\cong C_\lambda^*(S\rtimes S_{\n,\Lambda})$. Then
\setlength{\parindent}{0cm} \setlength{\parskip}{0cm}
\begin{enumerate}
\item[(i)] $K$ and $L$ are arithmetically equivalent, and $K(\m)^{\bar{\Gamma}}$ and $L(\n)^{\bar{\Lambda}}$ are Kronecker equivalent;
\item[(ii)] if the class fields $\Kmg$ and $L(\n)^{\bar{\Lambda}}$ of $K$ and $L$ are both Galois over $\Qz$, then $\#\clmg = \#\Cl_{\n}^{\bar{\Lambda}}$;
\item[(iii)] if $K$ or $L$ is Galois, then $K=L$; in particular, $K$ is Galois if and only if $L$ is Galois;
\item[(iv)] if $K$ or $L$ is Galois and both the class fields $\Kmg$ and $L(\n)^{\bar{\Lambda}}$ are Galois over $\Qz$, then
\begin{enumerate}
\item[(a)] $K=L$;
\item[(b)] $K(\m)^{\bar{\Gamma}}=L(\n)^{\bar{\Lambda}}$ (in any algebraically closed field containing both $K(\m)^{\bar{\Gamma}}$ and $L(\n)^{\bar{\Lambda}}$);
\item[(c)] $R^*\cdot(R_\n\cap R_{\m,\Gamma})=S^*\cdot (S_\m\cap S_{\n,\Lambda})$;
\item[(d)] $\clmg\cong\Cl_{\n}^{\bar{\Lambda}}$ (as abelian groups).
\end{enumerate}
\end{enumerate}
\etheooz
\setlength{\parindent}{0cm} \setlength{\parskip}{0cm}
We refer the reader to \S~\ref{s:Pre} for more precise definitions and more detailed explanations of our constructions. Here and in the sequel, we consider arithmetic equivalence and Kronecker equivalence over $\Qz$.
\setlength{\parindent}{0cm} \setlength{\parskip}{0.5cm}

We would like to point out that even in the case of trivial initial data, i.e., for ordinary $ax+b$-semigroups, our main theorem improves existing results. Indeed, in the present paper (more precisely in Theorem~\ref{thm:rootsofunity} and Corollary~\ref{cor:rootsofunity} below), we answer the natural question left open from \cite{Li14} whether it is possible to read off the number of roots of unity from our semigroup C*-algebras. This requires new techniques, including a detailed analysis of K-theoretic invariants, which is interesting on its own right. In addition, we are able to extract new information about the class fields naturally associated with our initial data (which in the case of trivial data is given by the Hilbert class fields). Indeed, by computing the torsion order of the $K_0$-class of the identity projection in quotients of $C^*_{\lambda}(R \rtimes R_{\m,\Gamma})$ by its minimal primitive ideals, we show that a certain set of prime numbers is an invariant of the C*-algebra $C^*_{\lambda}(R \rtimes R_{\m,\Gamma})$. This set differs from the Kronecker set of the class field $K(\m)^{\bar{\Gamma}}$ by only finitely many primes, so that we recover $K(\m)^{\bar{\Gamma}}$ up to Kronecker equivalence. This information was not available in \cite{Li14}. 

This paper is organized as follows. We begin with a brief discussion of general semigroup C*-algebras in \S~\ref{ss:SgpC} and then specialize to the case of semigroups of $ax+b$ type arising from actions of congruence monoids on rings of algebraic integers in \S~\ref{ss:conmon}. In \S~\ref{sec:cft} we explain how the number-theoretic data used to define a congruence monoid naturally gives rise to a class field (i.e., finite abelian extension), and in \S~\ref{ss:recnt} we show that how one can recover information about the congruence monoids from these class fields. In \S~\ref{s:pi}, we show that all the semigroup C*-algebras $C^*_{\lambda}(R \rtimes R_{\m,\Gamma})$ are purely infinite in a very strong sense, i.e., they absorb $\cO_{\infty}$ tensorially, $C^*_{\lambda}(R \rtimes R_{\m,\Gamma}) \cong \cO_{\infty} \otimes C^*_{\lambda}(R \rtimes R_{\m,\Gamma})$. This means that our semigroup C*-algebras fall into the scope of the classification programme for C*-algebras. More precisely, they belong to the class of C*-algebras classified by Kirchberg in \cite{Kir}. This motivates a detailed analysis of the K-groups of our semigroup C*-algebras, which we initiate in \S~\ref{ss:KSgpC}.
Furthermore, to illustrate our results, we discuss several concrete example classes along the way (see for instance \S~\ref{sss:K=Q}, \S~\ref{sec:realquadratic}, or Theorem~\ref{thm:Q}). We also present a first example of a semigroup whose left and right boundary quotients do not have isomorphic K-theory (see Remark~\ref{rem:lbq--rbq}), which is a completely new phenomenon that has not appeared before. This is an outgrowth of the general discussion of K-theory for boundary quotients of our semigroup C*-algebras in \S~\ref{ss:Kbq}. Finally, in \S~\ref{ss:RecCartan}, we prove stronger reconstruction results if in addition to the semigroup C*-algebras, we also keep track of canonical Cartan subalgebras. All these additional results show that our semigroup C*-algebras form an interesting class of examples from a C*-algebra perspective, giving access to many new examples of classifiable C*-algebras.

This project was initiated during a visit of the first-named author to Queen Mary University of London, and he would like to acknowledge this and thank the mathematics department there for its hospitality.

\section{Preliminaries}
\label{s:Pre}

\subsection{Semigroup C*-algebras}
\label{ss:SgpC}

Let $P$ be a left cancellative semigroup (with identity, say, for convenience), and consider the Hilbert space $\ell^2(P)$ with its canonical orthonormal basis $\{\delta_x :x\in P\}$. Since $P$ is left cancellative, each $p\in P$ gives rise to an isometry $\lambda(p):\ell^2(P)\to \ell^2(P)$ that is determined by $\lambda(p)(\delta_x)=\delta_{px}$ for $x\in P$. The \emph{left regular C*-algebra of $P$} is $C_\lambda^*(P):=C^*(\{\lambda(p):p\in P\})$. Let $I_l(P)$ be the inverse semigroup generated by the isometries $\lambda(P)$, and put $D_\lambda(P):=C^*(\{ss^* : s\in I_l(P)\})$. Each projection $ss^*$ for $s\in I_l(P)$ corresponds to a subset of $P$; such subsets are called constructible right ideals, and the set of constructible right ideals in $P$ is denoted by $\cJ_P$. 

Assume that $P$ embeds into a group $G$. Then $D_\lambda(P)$ coincides with the canonical diagonal sub-C*-algebra of $C_\lambda^*(P)$, namely $D_\lambda(P)=\ell^\infty(P)\cap C_\lambda^*(P)$ where we view $\ell^\infty(P)$ as a sub-C*-algebra of $\mathcal{L}(\ell^2(P))$ in the canonical way. 
Moreover, there is a canonical partial action of $G$ on $ D_{\lambda}(P)$, and $C^*_{\lambda}(P)$ can be written as the partial crossed product $C^*_{\lambda}(P) \cong D_{\lambda}(P) \rtimes_r G$. The C*-algebra $D_{\lambda}(P) \rtimes_r G$ can be identified with the reduced C*-algebra of the partial transformation groupoid $G\ltimes \Omega_P$ where $\Omega_P:=\Spec(D_{\lambda}(P))$, so one obtains a description of $C^*_{\lambda}(P)$ as a groupoid C*-algebra (see \cite{Li17} for details). 
This is useful for many purposes, for instance, in the case that $G\ltimes \Omega_P$ is topologically principal, so that $C(\Omega_P)$ is a Cartan subalgebra of $C_r^*(G\ltimes \Omega_P)$, we obtain that $D_\lambda(P)$ is a Cartan subalgebra of $C^*_{\lambda}(P)$.

We refer the reader to \cite{Li12,Li13,Li17} and \cite[Chapter~5]{CELY} for the general theory of semigroup C*-algebras.

\subsection{Congruence monoids and semigroup C*-algebras associated with their actions}
\label{ss:conmon}

We briefly review the construction from \cite[\S~3]{Bru}. Let $K$ be a number field with ring of algebraic integers $R$, and let $\cP_K$ denote the set of non-zero prime ideals of $R$. Each fractional ideal $\a$ of $K$ has a unique factorization $\a=\prod_{\p\in\cP_K}\p^{v_\p(\a)}$ where $v_\p(\a)\in\Zz$ is zero for all but finitely many $\p$. For $x\in K^\times$, we let $v_\p(x):=v_\p((x))$ where $(x)$ is the principal fractional ideal of $K$ generated by $x$. Given a modulus $\m=\m_\infty\m_0$ for $K$, let
\[
(R/\m)^*:=\bigg(\prod_{w\mid\m_\infty}\langle\pm 1\rangle\bigg)\times(R/\m_0)^*
\]
be the multiplicative group of residues modulo $\m$. Let 
\[
R_\m:=\{a\in R^\times :v_\p(a)=0 \text{ for all } \p\mid \m_0\}
\]
denote the multiplicative semigroup of non-zero algebraic integers that are relatively prime to $\m_0$; note that $R_\m$ only depends on the support of $\m_0$, that is, on $\supp(\m_0):=\{\p\in\cP_K : \p\mid \m_0\}$. For $a\in R_\m$, let
\[
[a]_\m:=((\sign(w(a)))_{w\mid\m_\infty},a+\m_0)\in(R/\m)^*.
\]
Then $R_\m\to (R/\m)^*$, $a\mapsto [a]_\m$, is a semigroup homomorphism. If $\Gamma$ is a subgroup of $(R/\m)^*$, then
\[
R_{\m,\Gamma}:=\{a\in R_\m :  [a]_\m\in\Gamma\}
\]
is a unital subsemigroup of $R^\times$, which is called a \emph{congruence monoid}. Note there is some freedom in choosing $\m_{\infty}$ and the infinite part of $\Gamma$ without changing the congruence monoid.

The monoid $R_{\m,\Gamma}$ acts on $R$ by multiplication, so one may form the semi-direct product $R\rtimes R_{\m,\Gamma}$. Let $K_{\m,\Gamma}:=\{x\in K_\m: [x]_\m\in\Gamma\}$; by \cite[Proposition~3.2]{Bru}, $K_{\m,\Gamma}=R_{\m,\Gamma}^{-1}R_{\m,\Gamma}$. By \cite[Proposition~3.3]{Bru}, the semigroup $R\rtimes R_{\m,\Gamma}$ is left Ore with group of left quotients equal to $G(R\rtimes R_{\m,\Gamma}) = (R_\m^{-1}R)\rtimes K_{\m,\Gamma}$, where $R_\m^{-1}R=\{a/b : a\in R, b\in R_\m\}\subseteq K^\times$ is the localization of $R$ at $R_\m$.

By \cite[Proposition~3.4]{Bru}, $R\rtimes R_{\m,\Gamma}$ satisfies the independence condition from \cite{Li12}, and the semilattice $\cJ_{R\rtimes R_{\m,\Gamma}}$ of constructible right ideals in $R\rtimes R_{\m,\Gamma}$ is given by
\[
\cJ_{R\rtimes R_{\m,\Gamma}}=\{(x+\a) \times (\a\cap \rmg) : x\in R, \a\in\cI_\m^+\}\cup\{\emptyset\}
\]
where $\cI_\m^+$ is the semigroup of integral ideals  relatively prime to $\m_0$. 

By \cite[equation~(3) and Proposition~5.4]{Bru}, we can identify $C^*_{\lambda}(R \rtimes M)$ with the reduced groupoid C*-algebra of the partial transformation groupoid $(Q \rtimes G) \ltimes \Omega$ in such a way that $D_\lambda(R\rtimes M)$ is carried onto $C(\Omega)$, where the space $\Omega$ and the partial action of $Q \rtimes G$ are described in \cite[\S~5.2]{Bru}. 
The groupoid $(Q \rtimes G) \ltimes \Omega$ is topologically principal by \cite[Proposition~6.3]{Bru}, so $D_\lambda(R\rtimes M)$ is a Cartan subalgebra of $C^*_{\lambda}(R \rtimes M)$ (cf. \S~\ref{ss:SgpC}).

Moreover, using the above description as a groupoid C*-algebra, the primitive ideal space of $C^*_{\lambda}(R \rtimes R_{\m,\Gamma})$ has been computed in \cite[\S~7]{Bru} as ${\rm Prim}(C^*_{\lambda}(R \rtimes R_{\m,\Gamma})) \cong 2^{\cP_K^{\mfm}}$, where $\cP_K^{\mfm}:=\cP_K\setminus\supp(\m_0)$, and $2^{\cP_K^{\mfm}}$ is given the power-cofinite topology. This homeomorphism is order-preserving, so it follows that the non-zero minimal primitive ideals of $C^*_{\lambda}(R \rtimes R_{\m,\Gamma})$ are in one-to-one correspondence with the primes in $\cP_K^{\mfm}$.

\subsection{Class fields associated with congruence monoids}\label{sec:cft}

We will need some standard results on ray class fields from the ideal-theoretic point of view. The reader may consult for instance \cite{MilCFT} for more details about class field theory. 

Let $K$ be a number field with ring of algebraic integers $R$, and let $\cP_K$ denote the set of non-zero prime ideals of $R$. Let $\m$ be a modulus for $K$, let $\cI_\m$ denote the group of fractional ideals of $K$ which are relatively prime to $\m_0$, and let $K_\m$ be the subgroup of $K^\times$ consisting elements that are  relatively prime to $\m_0$. Then the map $a\mapsto [a]_\m$ extends to a surjective group homomorphism $K_\m\to (R/\m)^*$ (see, for example, \cite[\S~2.2]{Bru}).
Let $i:K_\m\to \cI_\m$ be the canonical homomorphism given by $a\mapsto aR$. Let $K_{\m,1}:=\{x\in K_\m : [x]_\m=1\}$, so that $K_\m/K_{\m,1}\cong ( R/\m)^*$. The group $\Cl_\m(K):=\cI_\m/ i(K_{\m,1})$ is the \emph{ray class group modulo $\m$}. When it will not cause confusion, we will simply write $\Cl_\m$ rather than $\Cl_\m(K)$. The ray class group associated with the trivial modulus is the usual ideal class group of $K$, that is, $\Cl_{(1)}=\Cl$.

There is the following exact sequence relating $\Cl_\m$ to the ideal class group of $K$:
\bgl\label{gl:classgroups}
0\to R^*/R_{\m,1}^*\to K_\m/K_{\m,1}\to\Cl_\m\to\Cl\to 0
\egl
where $R_{\m,1}^*:=R^*\cap K_{\m,1}$ (see \cite[Chapter~V,~Theorem~1.7]{MilCFT}). 

Suppose $L/K$ is a finite extension of number fields. For a prime $\mathfrak{P}\in\cP_L$ lying over a prime $\p\in\cP_K$, we write $f_{L/K}(\mathfrak{P}\vert \p)$ for the inertia degree of $\mathfrak{P}$ over $\p$, and $e_{L/K}(\mathfrak{P}\vert \p)$ for the ramification index of $\mathfrak{P}$ over $\p$ (see \cite[Chapter~I,~\S8]{Neu} for the definitions). If $L$ is Galois, then $f_{L/K}(\mathfrak{P}\vert \p)$ does not depend on $\mathfrak{P}$, and we simply write $f_{L/K}(\p)$ instead of $f_{L/K}(\mathfrak{P}\vert \p)$.
When $K=\Qz$ and $\p=p\Zz$ for a rational prime $p$, we will often abuse notation and write $p$ instead of $p\Zz$.

Now suppose that $L/K$ is a finite Galois extension, and let $S$ denote the ring of algebraic integers in $L$. If a prime $\p\in\cP_K$ is unramified in $L$ and $\fP\in\cP_L$ with $\fP\mid\p$, then each element of the decomposition group 
\[
D_\fP(L/K)=\{\sigma\in\Gal(L/K) : \sigma(\fP)=\fP\}
\] 
defines an automorphism of the finite field $S/\fP$, and this gives a canonical identification of $D_\fP(L/K)$ with the Galois group $\Gal((S/\fP)/(R/\p))$ of the field extension $(S/\fP)/(R/\p)$. Since this latter group is cyclic, so is $D_\fP(L/K)$. The \emph{Frobenius automorphism corresponding to $\fP$} is the automorphism $(\fP,L/K)\in\Gal(L/K)$ such that under the identification $D_\fP(L/K)\cong\Gal((S/\fP)/(R/\p))$, $(\fP,L/K)$ is sent to the automorphism of $S/\fP$ that is determined by $x+\fP\mapsto x^{N(\p)}+\fP$. Here $N(\p) = \#(R/\p)$ is the norm of $\p$. It is known that the set 
\[
{\rm Fr}_\p(L/K):=\{(\fP,L/K) : \fP\in\cP_L \text{ with }\fP\mid\p\}
\]
is a conjugacy class in $\Gal(L/K)$. Indeed, $\Gal(L/K)$ acts transitively on the set of primes in $\cP_L$ that lie over $\p$, and if $\tau\in\Gal(L/K)$, then $(\tau(\fP),L/K)=\tau(\fP,L/K)\tau^{-1}$.

If $L/K$ is abelian, that is, if $L$ is a Galois extension of $K$ with abelian Galois group, then ${\rm Fr}_\p(L/K)$ consists of a single element, which we shall denote by $(\p,L/K)$; this automorphism is called the \emph{Frobenius automorphism corresponding to $\p$}.

For each modulus $\m$ of $K$, let $\Km$ denote the associated ray class field of $K$ (see \cite[Chapter~V,~\S~3]{MilCFT}). Then $\Km$ satisfies:
\setlength{\parindent}{0cm} \setlength{\parskip}{0cm}
\begin{enumerate}
\item[(a)] $\Km$ is a finite abelian extension of $K$,
\item[(b)] if a prime $\p\in\cP_K$ ramifies in $\Km$, then $\p\mid\m_0$, and if a real embedding $w$ of $K$ ramifies in $\Km$ --- meaning that there is a complex embedding of $\Km$ that extends $w$ --- then $w\mid\m_\infty$ (see \cite[Chapter~V,~Remark~3.8]{MilCFT});
\item[(c)] the Artin reciprocity map 
\[
\cI_\m\to \Gal(\Km/K),\, \a\mapsto\prod_{\p\mid\a}(\p,\Km/K)^{v_\p(\a)}
\] 
descends to a group isomorphism $\psi_{\Km/K}:\Cl_\m\overset{\cong}{\to}\Gal(\Km/K)$.
\end{enumerate}
Moreover, given moduli $\m$ and $\n$, we have $\Km\subseteq K(\n)$ if $\m \mid \n$.
\setlength{\parindent}{0cm} \setlength{\parskip}{0.5cm}

The ray class field $K(1)$ corresponding to the trivial modulus is called the \emph{Hilbert class field of $K$}. The Hilbert class field is the maximal everywhere unramified abelian extension of $K$. The Artin map gives an isomorphism $\Cl(K)\cong \Gal(K(1)/K)$, and we have the following commutative diagram with exact rows:
\bgl\label{gl:classgroup&galoisgroups}
\xymatrix{
0 \ar[r] & i(K_\m)/i(K_{\m,1})\ar[d]_{\cong}\ar[r] &  \Cl_\m\ar[d]_{\cong}^{\psi_{\Km/K}} \ar[r] & \Cl\ar[d]_{\cong}^{\psi_{K(1)/K}} \ar[r]& 0 \\
0 \ar[r] & \Gal(\Km/K(1))\ar[r] & \Gal(\Km/K)\ar[r] & \Gal(K(1)/K)\ar[r] & 0
}\egl

Using the isomorphism $K_\m/K_{\m,1}\cong (R/\m)^*$, we get a homomorphism $(R/\m)^*\cong K_\m/K_{\m,1}\to \Cl_\m$; by exactness of \eqref{gl:classgroups}, its kernel is $[R^*]_\m$, the image of the unit group $R^*$ in $(R/\m)^*$, and its range is $i(K_\m)/i(K_{\m,1})$. Therefore, using \eqref{gl:classgroup&galoisgroups}, we may identify $(R/\m)^*/[R^*]_\m$ with $\Gal(\Km/K(1))$. Thus, Galois theory gives an inclusion-reversing bijection between the set of subgroups of $(R/\m)^*/[R^*]_\m$ and the set of subfields of $K(\m)$ that contain $K(1)$.

Now we see that the arithmetic data $(\m,\Gamma)$ used to define the congruence monoid $R_{\m,\Gamma}$ also canonically determines a class field of $K$. Namely, let $\bar{\Gamma}$ be the image of $\Gamma$ under the composite $(R/\m)^* \onto (R/\m)^*/[R^*]_\m \cong \Gal(\Km/K(1))$, and consider the fixed field associated with $\bar{\Gamma}$, that is, the field $\Kmg$ consisting of elements in $\Km$ that are fixed by every element of $\bar{\Gamma}$. Observe that $\Kmg$ only depends on the image of $\Gamma$ in $(R/\m)^*/[R^*]_\m$ and that $\Kmg$ always contains $K(1)$.


Using the inclusion $ \Gal(\Km/K(1))\hookrightarrow \Gal(\Km/K)$, we may also view $\bar{\Gamma}$ as subgroup of $\Gal(\Km/K)$. The isomorphism $\psi_{\Km/K}:\Cl_\m\cong \Gal(\Km/K)$ takes the subgroup $i(K_{\m,\Gamma})/i(K_{\m,1})$ onto $\bar{\Gamma}$, so we have isomorphisms
\begin{equation}\label{eqn:galoisgroupisom}
\clmg:=\cI_\m/i(K_{\m,\Gamma})\cong \Cl_\m/(i(K_{\m,\Gamma})/i(K_{\m,1}))\cong\Gal(\Km/K)/\bar{\Gamma}\cong \Gal(\Kmg/K).
\end{equation}

\subsection{Reconstruction of initial data in number-theoretic context}
\label{ss:recnt}

We now give several results that will be used in the proofs of our reconstruction theorems. We are interested in the question how much of the congruence monoid $\rmg$ we can recover from $\Kmg$.

\bremark
The map $(\m,\Gamma)\mapsto R_{\m,\Gamma}$ is far from injective. However, from a number-theoretic perspective, it is natural to fix $\m$ and let $\Gamma$ vary; in terms of class fields, this corresponds to considering the intermediate extensions $K(1)\subseteq L\subseteq K(\m)$ for a fixed modulus $\m$.
\eremark
\setlength{\parindent}{0cm} \setlength{\parskip}{0cm}

We will continue using the notation from Section~\ref{sec:cft}. As we have seen, the number-theoretic data $(\m,\Gamma)$ used to define the congruence monoid $R_{\m,\Gamma}$ also defines a class field of $K$, namely the intermediate extension $K(1)\subseteq \Kmg\subseteq\Km$. From the discussion following \eqref{gl:classgroup&galoisgroups}, the extension $\Kmg$ only depends on the image of $\Gamma$ under the quotient map $(R/\m)^*\to (R/\m)^*/[R^*]_\m$, so one should not expect to be able to recover the monoid $R_{\m,\Gamma}$ from $\Kmg$. The following result shows exactly how much information is lost when passing from $R_{\m,\Gamma}$ to $\Kmg$.
\setlength{\parindent}{0cm} \setlength{\parskip}{0.5cm}

\bprop\label{prop:monoidfromclassfield}
Suppose that $\m$ and $\n$ are moduli for $K$ and that $\Gamma$ and $\Lambda$ are subgroups of $(R/\m)^*$ and $(R/\n)^*$, respectively. Then $K(\m)^{\bar{\Gamma}}=K(\n)^{\bar{\Lambda}}$ (equality in any algebraically closed field containing both $K(\m)^{\bar{\Gamma}}$ and $K(\n)^{\bar{\Lambda}}$) if and only if $R^*\cdot(R_\n\cap R_{\m,\Gamma})=R^*\cdot (R_\m\cap R_{\n,\Lambda})$. In particular, if $\supp(\m_0)=\supp(\n_0)$, then $K(\m)^{\bar{\Gamma}}=K(\n)^{\bar{\Lambda}}$ if and only if $R^*\cdot R_{\m,\Gamma}=R^*\cdot  R_{\n,\Lambda}$.
\eprop
\setlength{\parindent}{0cm} \setlength{\parskip}{0cm}

\bproof
It follows from \cite[Theorem~3.5.1]{Coh} that $K(\m)^{\bar{\Gamma}}=K(\n)^{\bar{\Lambda}}$ if and only if $\cI_\n\cap i(K_{\m,\Gamma})=\cI_\m\cap i(K_{\n,\Lambda})$. This is equivalent to having $\cI_\n^+\cap i(K_{\m,\Gamma})=\cI_\m^+\cap i(K_{\n,\Lambda})$, and since $\cI_n^+\cap i(K_{\m,\Gamma})=\cI_\n^+\cap i(R_{\m,\Gamma})$ and $\cI_\m^+\cap i(K_{\n,\Lambda})=\cI_\m^+\cap i(R_{\n,\Lambda})$, this is also equivalent to $\cI_\n^+\cap i(R_{\m,\Gamma})=\cI_\m^+\cap i(R_{\n,\Lambda})$. Since $a\in R^*\cdot(R_\n\cap R_{\m,\Gamma})$ if and only if $(a)\in \cI_\n^+\cap i(R_{\m,\Gamma})$, and $a\in R^*\cdot (R_\m\cap R_{\n,\Lambda})$ if and only if $(a)\in\cI_\m^+\cap i(R_{\n,\Lambda})$, we are done.
\eproof
\setlength{\parindent}{0cm} \setlength{\parskip}{0.5cm}

The following observation shows that, if we know $\supp(\m_0)$, then the class field $\Kmg$ remembers as much as possible about $R_{\m,\Gamma}$.

\blemma
We have $R^*\cdot R_{\m,\Gamma}=R_{\m,[R^*]_\m\cdot\Gamma}$ where $[R^*]_\m\cdot\Gamma$ is the subgroup of $(R/\m)^*$ generated by $\Gamma$ and $[R^*]_\m$, the image of the unit group.
\elemma
\setlength{\parindent}{0cm} \setlength{\parskip}{0cm}

\bproof
The containment ``$\subseteq$'' is obvious. Suppose that $a\in R_\m$ such that $[a]_\m\in [R^*]_\m\cdot\Gamma$. Then there exists $u\in R^*$ and $b\in R_{\m,\Gamma}$ such that $[a]_\m=[ub]_\m$. That is,
\begin{itemize}
\item $(\sign(w(a)))_{w\mid\m_\infty}=(\sign(w(u)))_{w\mid\m_\infty}(\sign(w(b)))_{w\mid\m_\infty}$, and 
\item $a+\m_0=ub+\m_0$.
\end{itemize}
By the second item above, there exists $x\in\m_0$ such that $a=ub+x=u(b+u^{-1}x)$. Now the first item implies that $(\sign(w(b+u^{-1}x)))_{w\mid\m_\infty}=(\sign(w(b)))_{w\mid\m_\infty}$, so we have $[b+u^{-1}x]_\m=[b]_\m$, that is, $b+u^{-1}x\in R_{\m,\Gamma}$.
\eproof
\setlength{\parindent}{0cm} \setlength{\parskip}{0.5cm}

We now turn to the natural question whether conversely $R_{\m,\Gamma}$ determines $\Kmg$.

\blemma\label{lem:splitsiff}
For each prime $\p$ not dividing $\m_0$, $\p$ is unramified in $\Kmg$. Let $f_\p^\Gamma$ denote the order of $[\p]$ in $\clmg$. Then $f_\p^\Gamma=f_{\Kmg/K}(\p)$, and $f_\p^\Gamma=1$ if and only if $\p$ splits completely in $\Kmg$.
\elemma
\setlength{\parindent}{0cm} \setlength{\parskip}{0cm}
\bproof
If $\p$ is  relatively prime to $\m_0$, then $\p$ is unramified in $\Km$ by property (b) of $K(\m)$ from \S~2.3, and thus $\p$ is also unramified in $\Kmg$.
\setlength{\parindent}{0cm} \setlength{\parskip}{0.5cm}

Under the isomorphism from \eqref{eqn:galoisgroupisom}, $[\p]$ is taken to  the Frobenius automorphism $(\p,\Kmg/K)$ whose order is precisely $f_{\Kmg/K}(\p)$. Since $\p$ splits in $\Kmg$ if and only if $(\p,\Kmg/K)=\id$, we are done.
\eproof
\setlength{\parindent}{0cm} \setlength{\parskip}{0.5cm}

\bprop\label{prop:isomofclassfields}
Let $\m$ and $\n$ be moduli for $K$, and let $\Gamma$ and $\Lambda$ be subgroups of $(R/\m)^*$ and $(R/\n)^*$, respectively. For $\p\in\cP_K^{\m\n}$, let $f_\p^\Gamma$ denote the order of $[\p]$ in $\clmg$ and let $f_\p^\Lambda$ be the order of $[\p]$ in $\clnl$. If $f_\p^\Gamma=f_\p^\Lambda$ for all but finitely many $\p$ in $\cP_K^{\m\n}$, then $\Kmg=\Knl$ (equality in any algebraic closure of $K$ that contains both $\Kmg$ and $\Knl$).
\eprop
\setlength{\parindent}{0cm} \setlength{\parskip}{0cm}

\bproof
Suppose that $f_\p^\Gamma=f_\p^\Lambda$ for all but finitely many $\p$ in $\cP_K^{\m\n}$. Then \Cref{lem:splitsiff} implies that for all but finitely many $\p$ in $\cP_K^{\m\n}$, $\p$ splits completely in $\Kmg$ if and only if $\p$ splits completely in $\Knl$. Hence, $\Kmg=\Knl$ by \cite[Chapter~V,~Theorem~3.25]{MilCFT}.
\eproof
\setlength{\parindent}{0cm} \setlength{\parskip}{0.5cm}

\bcor
Suppose $\m$ and $\n$ are moduli for $K$ and $\Gamma$ and $\Lambda$ are subgroups of $(R/\m)^*$ and $(R/\n)^*$, respectively. If $R_{\m,\Gamma}=R_{\n,\Lambda}$, then $\Kmg=\Knl$ (equality in any algebraically closed field containing both $K(\m)^{\bar{\Gamma}}$ and $K(\n)^{\bar{\Lambda}}$).
\ecor
\setlength{\parindent}{0cm} \setlength{\parskip}{0cm}

\bproof
If $R_{\m,\Gamma}=R_{\n,\Lambda}$, then $\Cl_\m^{\bar{\Gamma}}=\Cl_\n^{\bar{\Lambda}}$, so that for each $\p\in \cP_K^{\m\n}$, $f_\p^\Gamma=1$ if and only if $f_\p^\Lambda=1$. Hence, $\Kmg=\Knl$ by Proposition~\ref{prop:isomofclassfields}.
\eproof
\setlength{\parindent}{0cm} \setlength{\parskip}{0.5cm}

The class field $\Kz \defeq K(\mfm)^{\bar{\Gamma}}$ is always Galois over $K$, but we will need to know when $\Kz$ is Galois over $\Qz$. Let $\overline{\Qz}\subseteq \Cz$ denote the field of algebraic numbers.

\bprop\label{prop:congjugateclassfield}
Let $K$ be a number field. Suppose that $\m$ is a modulus for $K$ and that $\Gamma$ a group of residues modulo $\m$. Let $K\subseteq \overline{\Qz}$ be any embedding, so that we can view $\Kz$ as a subfield of $\overline{\Qz}$. If $\sigma\in\Gal(\overline{\Qz}/\Qz)$, then we denote by $\sigma(\m)$ the modulus for $\sigma(K)$ defined by $\sigma(\m)_0:=\sigma(\m_0)$ and $w\mid \sigma(\m)_\infty$ if and only if $w\circ\sigma\mid \m_\infty$.
Then $\sigma(\Kz)$ is the class field of $\sigma(K)$ corresponding to $(\sigma(\m),\sigma(\Gamma))$.
\eprop
\setlength{\parindent}{0cm} \setlength{\parskip}{0cm}

\bproof
Let $L=\sigma(K)$ and let $S$ denote the ring of algebraic integers in $L$. The isomorphism $R\to S$, $a\mapsto \sigma(a)$, defines an isomorphism $(R/\m)^*\cong (S/\sigma(\m))^*$ that we shall also denote by $\sigma$. Observe that $\sigma([R^*]_\m)=[S^*]_{\sigma(\m)}$, so that $\overline{\sigma(\Gamma)}=\sigma(\bar{\Gamma})$.
\setlength{\parindent}{0cm} \setlength{\parskip}{0.5cm}

We need to show that $\sigma(\Kz)$ and $\mathbb{L}:=L(\sigma(\m))^{\sigma(\bar{\Gamma})}$ are equal (as subfields of $\overline{\Qz}$). Now $\sigma$ induces an isomorphism $\cI_K^\m\cong \cI_L^{\sigma(\m)}$ where $\cI_K^\m$ denotes the group of fractional ideals of $K$ relatively prime to $\m_0$ and $\cI_L^{\sigma(\m)}$ denotes the group of fractional ideals of $L$ relatively prime to $\sigma(\m)_0$.
If $x\in K^\times$, then $x$ is relatively prime to $\m_0$ if and only if $\sigma(x)$ is relatively prime to $\sigma(\m)_0$, and $[x]_\m\in\Gamma$ if and only if $[\sigma(x)]_{\sigma(\m)}\in\sigma(\Gamma)$. So the isomorphism $K\cong L$ carries $K_{\m,\Gamma}$ onto $L_{\sigma(\m),\sigma(\Gamma)}$. Thus, we have an isomorphism $\Cl_\m^{\bar{\Gamma}}(K)\cong \Cl_{\sigma(\m)}^{\sigma(\bar{\Gamma})}(L)$ given by $[\a]\mapsto [\sigma(\a)]$. It follows that $\q\in\cP_L^{\sigma(\m)}$ splits completely in $\mathbb{L}$ if and only if $\sigma^{-1}(\q)$ splits completely in $\Kz$.

For $\a\in \cI_K^{\m}$, we have $\psi_{\sigma(\Kz)/\sigma(K)}(\sigma(\a))=\sigma\circ \psi_{\Kz/K}(\a)\circ\sigma^{-1}$ (see, for example, \cite[Chapter~X,~\S~1]{Lang}). Hence, the following diagram commutes:
\[\xymatrix{
\Cl_\m^{\bar{\Gamma}}(K)\ar[d]_{\cong}^{[\a]\mapsto [\sigma(\a)]}\ar[rr]_{\cong}^{\psi_{\Kz/K}} & \ \ \ &\ar[d]_{\cong}^{\tau\mapsto \sigma\circ\tau\circ\sigma^{-1}} \Gal(\Kz/K)\\
\Cl_{\sigma(\m)}^{\sigma(\bar{\Gamma})}(L)\ar[rr]_{\cong}^{\psi_{\sigma(\Kz)/L}}  & \ \ \ &\Gal(\sigma(\Kz)/L).
}\]
Thus Lemma~\ref{lem:splitsiff} implies that $\p\in\cP_K^\m$ splits completely in $\Kz$ if and only if $\sigma(\p)$ splits completely in $\sigma(\Kz)$.

All in all, we have that $\q\in\cP_L^{\sigma(\m)}$ splits completely in $\mathbb{L}$ if and only if $\q$ splits completely in $\sigma(\Kz)$. Hence,  $\mathbb{L}=\sigma(\Kz)$ by \cite[Chapter~V,~Theorem~3.25]{MilCFT}.
\eproof
\setlength{\parindent}{0cm} \setlength{\parskip}{0.5cm}

\bcor\label{cor:kzgalois}
Suppose that $K$ is a finite Galois extension of $\Qz$, $\m$ is a modulus for $K$, and $\Gamma$ is a group of residues modulo $\m$.
Then $\Kz$ is Galois over $\Qz$ if and only if $R^*\cdot (R_{\sigma(\m)}\cap R_{\m,\Gamma})=R^*\cdot (R_\m\cap R_{\sigma(\m),\sigma(\Gamma)})$ for all $\sigma\in\Gal(K/\Qz)$.
Hence, if $\sigma(\m)=\m$ and $\sigma(\Gamma)=\Gamma$ for all $\sigma\in \Gal(K/\Qz)$, then $\Kz$ is Galois over $\Qz$. 
\ecor
\setlength{\parindent}{0cm} \setlength{\parskip}{0cm}

\bproof
Since $K$ is Galois over $\Qz$, we may view $K$ as a subfield of $\overline{\Qz}$, and thus also view $\Kz$ as a subfield of $\overline{\Qz}$. The field $\Kz$ is Galois over $\Qz$ if and only if $\sigma(\Kz)=\Kz$ for all $\sigma\in\Gal(\overline{\Qz}/\Qz)$. Combining Proposition~\ref{prop:congjugateclassfield} and Proposition~\ref{prop:monoidfromclassfield}, we see that $\sigma(\Kz)=\Kz$ for all $\sigma\in\Gal(\overline{\Qz}/\Qz)$ if and only if $R^*\cdot (R_{\sigma(\m)}\cap R_{\m,\Gamma})=R^*\cdot (R_\m\cap R_{\sigma(\m),\sigma(\Gamma)})$ for every $\sigma\in\Gal(\overline{\Qz}/\Qz)$.
Since $K$ is Galois, each $\sigma\in\Gal(\overline{\Qz}/\Qz)$ maps $K$ onto itself and thus determines an element of $\Gal(K/\Qz)$. Moreover, every element of $\Gal(K/\Qz)$ is the restriction of an element in $\Gal(\overline{\Qz}/\Qz)$, so the above condition is equivalent to $R^*\cdot (R_{\sigma(\m)}\cap R_{\m,\Gamma})=R^*\cdot (R_\m\cap R_{\sigma(\m),\sigma(\Gamma)})$ for every $\sigma\in\Gal(K/\Qz)$.
\eproof
\setlength{\parindent}{0cm} \setlength{\parskip}{0.5cm}

\bex
Let $K$ be a finite Galois extension of $\Qz$ and let $l\in\Zz_{>0}$. Let $\m_\infty$ be either trivial or consist of all real embeddings of $K$, and let $\m_0=lR$. Then $\sigma(\m)=\m$ for all $\sigma\in\Gal(K/\Qz)$. 
The inclusion $\Zz\hookrightarrow R$ descends to an inclusion $(\Zz/l\Zz)^*\hookrightarrow (R/lR)^*$; if $\Gamma_0$ is the image of any subgroup of $(\Zz/l\Zz)^*$ under this inclusion, then $\Gamma:=\prod_{w\mid\m_\infty}\langle\pm 1\rangle\times\Gamma_0$ is a $\Gal(K/\Qz)$-invariant subgroup of $(R/\m)^*$. Now the corresponding class field $\Kz$ is Galois over $\Qz$ by Corollary~\ref{cor:kzgalois}.
\eex

\section{Pure infiniteness}
\label{s:pi}

{\bf Notational conventions:} To simplify notations, we will from now on and throughout this paper -- when it is safe to do so -- drop sub- and superscripts and write $M \defeq R_{\mfm,\Gamma}$, $Q = M^{-1} \cdot R$, $G = M^{-1} M$, $\mu \defeq \mu_{\mfm,\Gamma} \defeq \tor(M^*)$, $m \defeq \# \mu$, $C \defeq {\rm Cl}_{\mfm}^{\bar{\Gamma}} \ (\defeq {\rm Cl}_{\mfm}(K) / \bar{\Gamma})$, $\Kz \defeq K(\mfm)^{\bar{\Gamma}}$, and $f(\mfp) \defeq f_{\mfp}^{\Gamma}$.

In this section, we show that $C^*_{\lambda}(R \rtimes M)$ is strongly purely infinite.
\btheo
\label{thm:pi}
For every congruence monoid $M$ as in \S~\ref{ss:conmon}, we have $C^*_{\lambda}(R \rtimes M) \cong \cO_{\infty} \otimes C^*_{\lambda}(R \rtimes M)$.
\etheo
\setlength{\parindent}{0cm} \setlength{\parskip}{0cm}

\bproof
First of all, as explained in \S~\ref{ss:conmon}, $C^*_{\lambda}(R \rtimes M)$ is isomorphic to the reduced groupoid C*-algebra of the partial transformation groupoid $(Q \rtimes G) \ltimes \Omega$. Now the same proof as for \cite[Theorem~4.6]{Li17} -- with the following slight modifications -- shows that $(Q \rtimes G) \ltimes \Omega$ is purely infinite in the sense of \cite{Mat}: Replace the ideal $J$ defined as $\bigcap_{i=1}^n I_i$ in the proof of \cite[Theorem~4.6]{Li17} by $J \defeq (\bigcap_{i=1}^n I_i) \cap \mfm_0$ (where $\mfm_0$ is as in \S~\ref{ss:conmon}). Then follow the proof of \cite[Theorem~4.6]{Li17} to find an element $a \in (1+J) \setminus R^*$. The point is that because of our modified definition of $J$, we will always be able to find a suitable power $a^{\ve}$ of $a$ such that $a^{\ve}$ lies in $M$ (and $a^{\ve}$ still lies in $(1+J) \setminus R^*$ because this set is multiplicatively closed). With $a^{\ve}$ in place of $a$, the same proof as for \cite[Theorem~4.6]{Li17} shows that $(Q \rtimes G) \ltimes \Omega$ is purely infinite. As observed in \cite[\S~4.2]{Mat}, this implies that $C^*_{\lambda}(R \rtimes M) \cong C^*_r((Q \rtimes G) \ltimes \Omega)$ is purely infinite. As explained in \S~\ref{ss:conmon}, the primitive ideal space of $C^*_{\lambda}(R \rtimes M)$ is given by $2^{\cP_K^{\mfm}}$. And since $2^{\cP_K^{\mfm}}$ has a basis for its topology of compact-open subsets given by $U_F \defeq \menge{T \in 2^{\cP_K^{\mfm}}}{T \cap F = \emptyset}$, where $F$ runs through all finite subsets of $\cP_K^{\mfm}$, and because $C^*_{\lambda}(R \rtimes M)$ is separable and purely infinite, it follows from \cite[Proposition~2.11]{PR} that $C^*_{\lambda}(R \rtimes M)$ has the ideal property from \cite{PR}. Hence \cite[Proposition~2.14]{PR} implies that $C^*_{\lambda}(R \rtimes M)$ is strongly purely infinite. Finally, as $C^*_{\lambda}(R \rtimes M)$ is separable, nuclear and unital, \cite[Theorem~8.6]{KR} implies that $C^*_{\lambda}(R \rtimes M) \cong \cO_{\infty} \otimes C^*_{\lambda}(R \rtimes M)$, as desired.
\eproof
\setlength{\parindent}{0cm} \setlength{\parskip}{0.5cm}

\section{K-theory}\label{sec:ktheory}

\subsection{K-theory for our semigroup C*-algebras}
\label{ss:KSgpC}

First of all, let us compute K-theory for our semigroup C*-algebras. For each class $\k\in C$, choose an integral ideal $\a_\k\in \k$; for $\k = [R]$, take $\a_\k=R$. Let $\iota_\k$ denote the homomorphism $C^*(\a_\k\rtimes M^*)\to C^*_{\lambda}(R \rtimes M)$ determined by $u_x\mapsto \lambda(x) e_{\a_\k}$ where, for $x\in\a_\k\rtimes M^*$, $u_x$ denotes the corresponding unitary in $C^*(\a_\k\rtimes M^*)$ and $e_{\a}$ is the projection in $C^*_{\lambda}(R \rtimes M)$ corresponding to the constructible right ideal $\a\times (\a\cap \rmg)$ of $R\rtimes M$.

\btheo\label{thm:ktheoryformula}
The homomorphisms $\iota_\k$ induce an isomorphism 
\[
\sum_{\k\in C}(\iota_\k)_*:\bigoplus_{\k\in C}K_*(C^*(\a_\k\rtimes M^*))\cong K_*(C^*_{\lambda}(R \rtimes M)).
\]
\etheo
\setlength{\parindent}{0cm} \setlength{\parskip}{0cm}

Here $K_*$ denotes $K_0\oplus K_1$ as a $\Zz/2\Zz$-graded abelian group.
\bproof
The group $G(R \rtimes M) = Q \rtimes G$ is solvable, hence amenable, so $Q \rtimes G$ satisfies the Baum-Connes conjecture with arbitrary coefficients by \cite{HK}. Since $R \rtimes M \subseteq Q \rtimes G$ is left Ore by \cite[Proposition~3.2]{Bru}, $R \rtimes M \subseteq Q \rtimes G$ satisfies the left Toeplitz condition; by \cite[Proposition~3.4]{Bru}, $\cJ_{R \rtimes M}$ is independent, so $\cJ_{R \rtimes M \subseteq Q \rtimes G}$ is independent by \cite[Lemma~4.2]{Li13}. Thus, we can apply \cite[Corollary~3.14]{CEL2} (take for the group $G$ in \cite{CEL2} our group $Q \rtimes G$ and for $I$ in \cite{CEL2} the set $\cJ_{R \rtimes M \subseteq Q \rtimes G}^\times$).
\eproof
\setlength{\parindent}{0cm} \setlength{\parskip}{0.5cm}

\bremark
There is a parallel between the K-theory formula for $C^*_{\lambda}(R \rtimes M)$ given by Theorem~\ref{thm:ktheoryformula} and the parameterization of the low temperature KMS states on $C^*_{\lambda}(R \rtimes M)$ with respect to the canonical time evolution coming from the norm map on $K$. 
Indeed, \cite[Theorem~3.2(iii)]{Bru2} implies that for each fixed $\beta>2$, the simplex of KMS$_\beta$ states on $C^*_{\lambda}(R \rtimes M)$ is isomorphic to the simplex of tracial states on the C*-algebra $\bigoplus_{\k\in C}C^*(\a_\k\rtimes M^*)$.
\setlength{\parindent}{0.5cm} \setlength{\parskip}{0cm}

This connection has been observed in the case of the full $ax+b$-semigroup $R\rtimes R^\times$ (see the discussion following \cite[Theorem~6.6.1]{CELY}).
\eremark
\setlength{\parindent}{0cm} \setlength{\parskip}{0.5cm}

\subsubsection{The case of the rational number field}
\label{sss:K=Q}

Consider the case $K=\Qz$. Then the group C*-algebras appearing in the K-theory formula given by Theorem~\ref{thm:ktheoryformula} are all isomorphic to either $\Zz\rtimes\langle\pm 1\rangle$ or $\Zz$ according to whether or not $-1\in\Zz_{\m,\Gamma}$. Hence,
\begin{itemize}
\item $K_{\bullet}(C^*_{\lambda}(\Zz\rtimes\Zz_{\m,\Gamma}))\cong 
\begin{cases}
\Zz^{3\cdot\#C} &\falls \bullet=0,\\
\{0\}  &\falls \bullet =1
\end{cases}$\quad when $-1\in\Zz_{\m,\Gamma}$;
\item $K_{\bullet}(C^*_{\lambda}(\Zz\rtimes\Zz_{\m,\Gamma}))\cong \begin{cases}
\Zz^{\#C} & \falls \bullet=0,\\
\Zz^{\#C} & \falls \bullet =1
\end{cases}$ \quad when $-1\notin\Zz_{\m,\Gamma}$.
\end{itemize}

\subsubsection{On the summands in the K-theory formula for our semigroup C*-algebras}
\label{ss:Ksummands}

Let us now compare the summands in the K-theory formula in Theorem~\ref{thm:ktheoryformula}.
\bprop
\label{PROP:I=R}
For every non-zero ideal $\mfa$ of $R$, we have that $K_*(C^*(\mfa \rtimes M^*))$ and $K_*(C^*(R \rtimes M^*))$ are isomorphic up to inverting $m$, i.e., 
$$
  \Zz[\tfrac{1}{m}] \otimes K_*(C^*(\mfa \rtimes M^*)) \cong  \Zz[\tfrac{1}{m}] \otimes K_*(C^*(R \rtimes M^*)).
$$
\eprop

For the proof, we recall the main result from \cite{LaLu}, as it is stated for $\mfa = R$ and even $m$ in \cite{LiLu}.

Let $\mfa$ be a non-zero ideal of $R$, $\iota: \mfa \to \mfa \rtimes \mu$ the canonical inclusion, and denote by $\iota_*$ the homomorphism $K_{\bullet}(C^*(\mfa)) \to K_{\bullet}(C^*(\mfa \rtimes \mu))$ induced by $\iota$ on $K_{\bullet}$ ($\bullet = 0,1$). Moreover, given a finite subgroup $F$ of $\mfa \rtimes \mu$, consider the canonical projection $F \onto \gekl{e}$ from $F$ onto the trivial group. This projection induces a homomorphism $C^*(F) \to \Cz$ of group C*-algebras, hence a homomorphism on $K_0$, $K_0(C^*(F)) \to K_0(\Cz)$. Let us denote the kernel of this homomorphism by $\widetilde{R}_\Cz(F)$. The canonical inclusion $F \to \mfa \rtimes \mu$ induces a homomorphism $\iota_F: C^*(F) \to C^*(\mfa \rtimes \mu)$, hence a homomorphism $K_0(C^*(F)) \to K_0(C^*(\mfa \rtimes \mu))$. Restricting this homomorphism to $\widetilde{R}_\Cz(F)$, we obtain $(\iota_F)_*: \widetilde{R}_\Cz(F) \to K_0(C^*(\mfa \rtimes \mu))$. Here are the main results from \cite{LaLu}:

\btheo[Langer-L\"uck]
\label{LL}
With the notations from above, we have
\setlength{\parindent}{0cm} \setlength{\parskip}{0cm}

\begin{itemize}
\item $K_0(C^*(\mfa \rtimes \mu))$ is finitely generated and torsion-free.

\begin{itemize}
\item 
$\rk_{\Zz} \rukl{\img(\iota_*)} = \rk_{\Zz} \rukl{\rukl{K_0(C^*(\mfa))}^{\mu}}$, and if $\cM$ denotes the set of conjugacy classes of maximal finite subgroups of $\mfa \rtimes \mu$, then $\sum_{(F) \in \cM} (\iota_F)_*: \bigoplus_{(F) \in \cM} \widetilde{R}_\Cz(F) \to K_0(C^*(\mfa \rtimes \mu))$ is injective.

\item
We have $\img(\iota_*) \cap \rukl{\sum_{(F) \in \cM} \img((\iota_F)_*)} = \gekl{0}$, and $\iota_* \oplus \rukl{\sum_{(F) \in \cM} (\iota_F)_*}$ is surjective after inverting $m$.
\end{itemize}

\item
$K_1(C^*(\mfa \rtimes \mu))$ is finitely generated and torsion-free.

\begin{itemize}
\item The map $K_1(C^*(\mfa))^{\mu} \to K_1(C^*(\mfa \rtimes \mu))$ induced by $\iota_*$ is an isomorphism after inverting $m$.
\end{itemize}

\end{itemize}
\etheo
\setlength{\parindent}{0cm} \setlength{\parskip}{0.5cm}

As above, let $\mfa$ be a non-zero ideal of $R$. We obtain canonical inclusions $\mfa \into R$ and $\mfa \rtimes \mu \into R \rtimes \mu$, both of which are denoted by $i$. Let $K_{\rm fin}^\mfa \subseteq K_0(C^*(\mfa \rtimes \mu))$ be defined by
$$
  K_{\rm fin}^\mfa \defeq \Big(\sum_{(F) \in \cM} \img((\iota_F)_*)\Big).
$$
Moreover, for $\bullet = 0,1$, let $K_{\bullet, \rm inf}^\mfa \subseteq K_{\bullet}(C^*(\mfa \rtimes \mu))$ be given by $K_{\bullet, \rm inf}^\mfa \defeq \img(\iota_*)$. Here we use the same notation as in the Theorem~\ref{LL}. Moreover, let $N(\mfa) \defeq \# R/\mfa$.

\blemma
\label{LEM:rkKinf}
For every non-zero ideal $\mfa$ of $R$, we have $\rk_{\Zz}(K_{\bullet, \rm inf}^\mfa) = \rk_{\Zz}(K_{\bullet, \rm inf}^R)$ for $\bullet = 0,1$.
\elemma
\setlength{\parindent}{0cm} \setlength{\parskip}{0cm}

\bproof
Clearly, $i_*: \: K_{\bullet}(C^*(\mfa)) \to K_{\bullet}(C^*(R))$ is injective and $\mu$-equivariant, so that it induces an embedding $\rukl{K_{\bullet}(C^*(\mfa))}^{\mu} \to \rukl{K_{\bullet}(C^*(R))}^{\mu}$. By Theorem~\ref{LL}, this shows $\rk_{\Zz}(K_{\bullet, \rm inf}^\mfa) \leq \rk_{\Zz}(K_{\bullet, \rm inf}^R)$. For the reverse inequality, take $x \in \rukl{K_{\bullet}(C^*(R))}^{\mu}$. Recall that $N(\mfa) = \# R/\mfa$. Set $n \defeq [K:\Qz] = \rk_{\Zz}(R)$. Then $N(\mfa)^n \cdot x$ lies in $\img(i_*)$, say $N(\mfa)^n \cdot x = i_*(y)$. Since $i_*$ is injective and $\mu$-equivariant, $y$ must be fixed by $\mu$. Hence the image of $\rukl{K_{\bullet}(C^*(\mfa))}^{\mu} \to \rukl{K_{\bullet}(C^*(R))}^{\mu}$ has finite index. This shows $\rk_{\Zz}(K_{\bullet, \rm inf}^\mfa) \geq \rk_{\Zz}(K_{\bullet, \rm inf}^R)$.
\eproof
\setlength{\parindent}{0cm} \setlength{\parskip}{0.5cm}

It is straighforward to check that $i_*(K_{\bullet, \rm inf}^\mfa) \subseteq K_{\bullet, \rm inf}^R$. Since $i_*: \: K_{\bullet}(C^*(\mfa)) \to K_{\bullet}(C^*(R))$ is surjective after inverting $N(\mfa)$, Lemma~\ref{LEM:rkKinf} has the following immediate consequence:
\bcor
\label{COR:Kinf}
$i_* \vert_{K_{\bullet, \rm inf}^\mfa}: \: K_{\bullet, \rm inf}^\mfa \to K_{\bullet, \rm inf}^R$ is an isomorphism after inverting $N(\mfa)$.
\ecor

Similarly, we have that $i_*(K_{\rm fin}^\mfa) \subseteq K_{\rm fin}^R$. Our goal is to show that for particular choices of $\mfa$, $i_* \vert_{K_{\rm fin}^\mfa}: K_{\rm fin}^\mfa \to K_{\rm fin}^R$ is an isomorphism.

Let $\pi: \: \mfa \rtimes \mu \onto \mu$ be the canonical projection.
\blemma
For every finite subgroup $F \subseteq \mfa \rtimes \mu$, $\pi \vert_F: \: F \to \mu$ is injective.
\elemma
\setlength{\parindent}{0cm} \setlength{\parskip}{0cm}

\bproof
Take $x \in \ker(\pi \vert_F)$. Then $x \in \ker(\pi) = \mfa$ and $x \in F$. However, every non-zero element in $\mfa$ has infinite order as $\mfa$ is torsion-free. Hence $x$ must be trivial, so that $\pi \vert_F$ is indeed injective.
\eproof
\setlength{\parindent}{0cm} \setlength{\parskip}{0.5cm}

Let $\zeta$ be a root of unity such that $\mu = \spkl{\zeta}$.
\bcor
Every finite subgroup $F$ is of the form $F = \spkl{(r,\zeta^i)}$ for some natural number $0 \leq i \leq m-1$ and $r \in \mfa$.
\ecor

\blemma
\label{LEM:conj}
Two elements $(r,\zeta^i)$ and $(s,\zeta^j)$ are conjugate in $\mfa \rtimes \mu$ if and only if $i=j$ and there exist $\xi \in \mu$ and $t \in \mfa$ such that $s = \xi(r + (1 - \zeta^i)t)$.
\elemma
\setlength{\parindent}{0cm} \setlength{\parskip}{0cm}

\bproof
$(r,\zeta^i)$ and $(s,\zeta^j)$ are conjugate if and only if there is $(b,a) \in \mfa \rtimes \mu$ such that
$$
  (s,\zeta^j) = (b,a)(r,\zeta^i)(b,a)^{-1} = (ar + (1 - \zeta^i)b,\zeta^i),
$$
which holds if and only if $i = j$ and $s = ar + (1 - \zeta^i)b$. Set $\xi \defeq a$ and $t \defeq a^{-1}b$.
\eproof
\setlength{\parindent}{0cm} \setlength{\parskip}{0.5cm}

Now take $\mfa$ such that, for every $1 \neq \xi \in \mu$, $\mfa$ and $(1 - \xi)$ are relatively prime (as ideals of $R$). Let $\cM^\mfa$ be the set of conjugacy classes of maximal finite subgroups of $\mfa \rtimes \mu$.
\blemma
\label{LEM:cMI=cMR}
For every such $\mfa$, $\cM^\mfa \to \cM^R, \, (F) \ma (F)$ is a bijection.
\elemma
\setlength{\parindent}{0cm} \setlength{\parskip}{0cm}

\bproof
We write $\sim_\mfa$ for conjugacy in $\mfa \rtimes \mu$.
\setlength{\parindent}{0cm} \setlength{\parskip}{0.5cm}

We first show that for a finite subgroup $F' = \spkl{(r,\xi)}$ of $R \rtimes \mu$, there exists a finite subgroup $\ti{F}$ of $\mfa \rtimes \mu$ such that $F \sim_R \ti{F}$. The point is that since $\mfa$ and $(1 - \xi)$ are relatively prime, we have $R = \mfa + (1 - \xi)$, so that there exists $t \in R$ and $s \in \mfa$ such that $r = s + (1 - \xi)t$. Now set $\ti{F} \defeq \spkl{(s,\xi)}$. Then $\ti{F} \subseteq \mfa \rtimes \mu$, and by Lemma~\ref{LEM:conj}, we have $F' \sim_R \ti{F}$.

Secondly, we show that given $r, s \in \mfa$ and $1 \neq \xi \in \mu$ with $(r,\xi) \sim_R (s,\xi)$, say $(s,\xi) = (b,a)(r,\xi)(b,a)^{-1}$, we must have $(b,a) \in \mfa \rtimes \mu$, i.e., $(r,\xi) \sim_\mfa (s,\xi)$. Namely, by the same computation as in Lemma~\ref{LEM:conj}, $(s,\xi) = (b,a)(r,\xi)(b,a)^{-1}$ implies that $s = ar + (1 - \xi)b$. But then $(1 - \xi) b = s - a r \in \mfa$, so that $b \in ((1 - \xi)^{-1} \cdot \mfa) \cap R = \mfa$ since $\mfa$ and $(1 - \xi)$ are relatively prime.

With these two observations, we show that the map $\cM^\mfa \to \cM^R, \, (F) \ma (F)$ is well-defined: Assume that $\gekl{e} \neq F \subseteq \mfa \rtimes \mu$ and $F' \subseteq R \rtimes \mu$ are finite subgroups with $F \subseteq F'$. By our first observation, there exists $\ti{F} \subseteq \mfa \rtimes \mu$ such that $F \sim_R \ti{F}$, say $F = (b,a) \ti{F} (b,a)^{-1}$. Suppose $F = \spkl{(r,\xi)}$, where $\xi \neq 1$ as $F \neq \gekl{e}$. As $F \subseteq F'$, we must have $(b,a)^{-1}(r,\xi)(b,a) \in \mfa \rtimes \mu$. By our second observation, this implies $(b,a) \in \mfa \rtimes \mu$, so that $F' \subseteq (b,a)(\mfa \rtimes \mu)(b,a)^{-1} = \mfa \rtimes \mu$. This shows that maximal finite subgroups of $\mfa \rtimes \mu$ are still maximal in $R \rtimes \mu$. So $\cM^\mfa \to \cM^R, \, (F) \ma (F)$ is well-defined.

By our first observation, $\cM^\mfa \to \cM^R, \, (F) \ma (F)$ is surjective. To see injectivity, let $F$ and $F'$ be maximal finite subgroups of $\mfa \rtimes \mu$, and suppose that $F \sim_R F'$, say $\gamma F \gamma^{-1} = F'$ for some $\gamma \in R \rtimes \mu$. Suppose $F = \spkl{(r,\xi)}$. Let $(s,\xi) = \gamma (r,\xi) \gamma^{-1}$. Then $F' = \spkl{(s,\xi)}$. In particular, $(s,\xi) \in \mfa \rtimes \mu$, and by our second observation, we must have $(r,\xi) \sim_\mfa (s,\xi)$, i.e., $F \sim_\mfa F'$.
\eproof
\setlength{\parindent}{0cm} \setlength{\parskip}{0.5cm}

\bcor
\label{COR:Kfin}
For every $\mfa$ as in Lemma~\ref{LEM:cMI=cMR}, $i_* \vert_{K_{\rm fin}^\mfa}: \: K_{\rm fin}^\mfa \to K_{\rm fin}^R$ is an isomorphism.
\ecor

We are now ready for
\setlength{\parindent}{0cm} \setlength{\parskip}{0cm}

\bproof[Proof of Proposition~\ref{PROP:I=R}]
By Corollary~\ref{COR:Kinf} and \ref{COR:Kfin}, for very ideal $\mfa$ as in Lemma~\ref{LEM:cMI=cMR}, we have in $K_0$ that
$$
  i_* \vert_{K_{0, \rm inf}^\mfa + K_{\rm fin}^\mfa}: \: K_{0, \rm inf}^\mfa + K_{\rm fin}^\mfa \to K_{0, \rm inf}^R + K_{\rm fin}^R
$$
is an isomorphism after inverting $N(\mfa)$. By Theorem~\ref{LL}, $i_*: \: K_0(C^*(\mfa \rtimes \mu)) \to K_0(C^*(R \rtimes \mu))$ is an isomorphism after inverting $m \cdot N(\mfa)$.

In $K_1$, consider the commutative diagram
$$
\xymatrix@C=1mm{
  K_1(C^*(\mfa)) \ar[d] \ar[rr]_{i_*} & \ \ \ & K_1(C^*(R)) \ar[d]\\
  K_1(C^*(\mfa \rtimes \mu)) \ar[rr]^{i_*} & \ \ \ & K_1(C^*(R \rtimes \mu))
}
$$
The upper horizontal arrow is bijective after inverting $N(\mfa)$, and the vertical arrows are surjective after inverting $m$, so that the lower horizontal arrow must be surjective after inverting $m \cdot N(\mfa)$. This together with Lemma~\ref{LEM:rkKinf} implies that the lower horizontal arrow must be an isomorphism after inverting $m \cdot N(\mfa)$.

Hence, by applying the Pimsner-Voiculescu sequence, we see that $i_*: \: K_*(C^*(\mfa \rtimes M^*)) \to K_*(C^*(R \rtimes M^*))$ is an isomorphism after inverting $m \cdot N(\mfa)$. In particular, $\rk_{\Zz}(K_*(C^*(\mfa \rtimes M^*))) = \rk_{\Zz}(K_*(C^*(R \rtimes M^*)))$. Now given an arbitrary non-zero ideal $\mfa$ of $R$, we can choose another ideal $\mfb$ of $R$ in the same ideal class of $\mfa$ such that $\mfb$ is relatively prime to $(1 - \xi)$ for all $1 \neq \xi \in \mu$, and such that $N(\mfb)$ is relatively prime to the order of every torsion element in $K_*(C^*(\mfa \rtimes M^*))$ and $K_*(C^*(R \rtimes M^*))$. Then $i_*$ must be already injective after inverting $m$, so that $i_*$ induces an isomorphism on the torsion parts after inverting $m$. Since we already know $\rk_{\Zz}(K_*(C^*(\mfa \rtimes M^*))) = \rk_{\Zz}(K_*(C^*(R \rtimes M^*)))$, we are done.
\eproof

\bremark
Already in Theorem~\ref{LL}, we have to invert $m$. This seems to be difficult to avoid.
\eremark

\bcor
\label{COR:rk=hrk}
We have $\rk(\Qz \otimes K_0(C^*_{\lambda}(R \rtimes M))) = (\# C) \cdot \rk(\Qz \otimes K_0(C^*_{\lambda}(R \rtimes M^*)))$. In particular, $\infty > \rk(\Qz \otimes K_0(C^*_{\lambda}(R \rtimes M))) \geq \#C$.
\ecor
\setlength{\parindent}{0cm} \setlength{\parskip}{0cm}

\bproof
The first claim follows immediately from the K-theory formula in Theorem~\ref{thm:ktheoryformula} and Proposition~\ref{PROP:I=R}. The second claim follows from the first once we know that $\infty > \rk(\Qz \otimes K_0(C^*_{\lambda}(R \rtimes M^*))) > 0$. The latter holds because the canonical maps $C^*(M^*) \to C^*(R \rtimes M^*) \to C^*(M^*)$ compose to the identity, so that $K_0(C^*(M^*)) \to K_0(C^*(R \rtimes M^*))$ is injective, and $K_0(C^*(M^*))$ is always free abelian of finite (but strictly positive) rank. This shows $\rk(\Qz \otimes K_0(C^*_{\lambda}(R \rtimes M^*))) > 0$. Moreover, we have $\rk(\Qz \otimes K_0(C^*_{\lambda}(R \rtimes M^*))) < \infty$ because of Theorem~\ref{LL} and the Pimsner-Voiculescu exact sequence.
\eproof
\setlength{\parindent}{0cm} \setlength{\parskip}{0.5cm}

\subsubsection{The case of real quadratic fields and totally positive elements}\label{sec:realquadratic}

Now consider a real quadratic field $K$ with ring of algebraic integers $R$. Suppose $K=\Qz(\sqrt{d})\subseteq\Rz$ where $d>1$ is a square-free natural number. Let $\mfm = \mfm_{\infty}$ be given by all infinite places of $K$, i.e., the real embeddings determined by $\sqrt{d} \ma \sqrt{d}$ and $\sqrt{d} \ma - \sqrt{d}$. Let $\Gamma \subseteq (R/\mfm)^*$ be the trivial subgroup. Then the congruence monoid $M = R_{\mfm,\Gamma}$ is given by $R_+\reg$, the set of (non-zero) totally positive elements in $R$. Our goal is to explicitly compute K-theory for $C^*_{\lambda}(R \rtimes R_+\reg)$.

The following result is a special case of the analysis from \cite{Is}, but we give a slightly different presentation here.  Let $\epsilon$ be the generator of $R_+^*$ with $\epsilon>1$. Then $\epsilon=\frac{t+u\sqrt{D}}{2}$ where $D$ is the discriminant of $R$ and $(t,u)$ is the smallest positive solution to the Pell equation $x^2-Dy^2=4$ (see \cite[Chapter~I,~\S~7,~Exercise~1]{Neu}.)

Let $\a$ be a fractional ideal of $K$, and denote by $\beta_\epsilon$ the automorphism of $C^*(\a)$ determined by $\beta_{\epsilon}(u_x)=u_{\epsilon x}$ where $u_x$ denotes the unitary in $C^*(\a)$ corresponding to $x\in\a$.

\bprop\label{prop:ktheoryrealquad}
 The induction map 
 \[
 \ind_{\a}^{\a\rtimes\langle \epsilon\rangle  }: K_0(C^*(\a))\to K_0(C^*(\a\rtimes\langle \epsilon\rangle  ))
 \]
is an isomorphism of $K_0$-groups, and there is an isomorphism 
\[
K_1(C^*(\a\rtimes\langle \epsilon\rangle  ))\cong K_0(C^*(\a))\oplus \a/(1-\epsilon)\a
\] 
of abelian groups. Indeed, the following sequence is exact:
\[
0\longrightarrow K_1(C^*(\a)) \overset{\id-(\beta_\epsilon)_*}{\longrightarrow} K_1(C^*(\a)) \overset{\ind_{\a}^{\a\rtimes\langle \epsilon\rangle  }}{\longrightarrow} K_1(C^*(\a\rtimes\langle \epsilon\rangle )) \overset{\partial}{\longrightarrow}  K_0(C^*(\a))\longrightarrow 0
\]
where $\partial$ is the boundary map from the Pimsner-Voiculescu exact sequence for $\beta_\epsilon:\Zz\acts C^*(\a)$.
\eprop
\setlength{\parindent}{0cm} \setlength{\parskip}{0cm}

\bproof
Let $w_1,w_2$ be a $\Zz$-basis for $\a$. Then $K_1(C^*(\a))\cong \Zz^2$ has a $\Zz$-basis given by $\{[u_{w_1}]_1,[u_{w_2}]_1\}$, and $K_0(C^*(\a))\cong \Zz^2$ has a $\Zz$-basis given by $\{[1]_0, [u_{w_1}]_1\times [u_{w_2}]_1\}$, were ``$\times$" denotes the product in K-theory (see \cite[Chapter~4.7]{HR}). Let $\gamma\in \textup{SL}_2(\Zz)$ be the matrix for $\beta_\epsilon$ with respect to the basis $w_1,w_2$. Then $(\beta_\epsilon)_*$ is simply given by applying the matrix $\gamma$ on $K_1$, and is trivial on $K_0$ since $(\beta_w)_*([1]_0)=[1]_0$, and $(\beta_\epsilon)_*([u_{w_1}]_1\times [u_{w_2}]_1)=\det(\gamma)([u_{w_1}]_1\times [u_{w_2}]_1)=[u_{w_1}]_1\times [u_{w_2}]_1$.
Now, $\textup{trace}(\gamma)=t$, and $t^2=4+Du^2$ with $t,u>0$, so we must have $t>2$.
Since $\det(\id-\gamma)=2-\textup{trace}(\gamma)$ is non-zero, we see that $\id-(\beta_\epsilon)_*$ is injective on $K_1(C^*(\a))$, so the Pimsner-Voiculescu exact sequence implies the results after noting that $K_1(C^*(\a))/\im(\id-(\beta_\epsilon)_*)\cong \a/(1-\epsilon)\a$ as abelian groups.
\eproof
\setlength{\parindent}{0cm} \setlength{\parskip}{0.5cm}

\btheo
\label{thm:KCrealquad}
Let $K$ be a real quadratic field with ring of algebraic integers $R$. Then
\[
K_0(C^*(R\rtimes R_+^\times))\cong\Zz^{2h_K^+} \quad\text{ and }\quad K_1(C^*(R\rtimes R_+^\times))\cong\Zz^{2h_K^+}\oplus (R/(1-\epsilon)R)^{h_K^+}
\]
where $h_K^+$ is the narrow class number of $K$.
\etheo
\setlength{\parindent}{0cm} \setlength{\parskip}{0cm}

\bproof
This follows from \Cref{thm:ktheoryformula}, Proposition~\ref{PROP:I=R} and \Cref{prop:ktheoryrealquad}.
\eproof
\setlength{\parindent}{0cm} \setlength{\parskip}{0.5cm}

Let us present an immediate consequence. Suppose $K$ and $K'$ are real quadratic fields with rings of algebraic integers $R$ and $R'$, and write $K=\Qz(\sqrt{d})$ and $K'=\Qz(\sqrt{d'})$ where $d$ and $d'$ are (uniquely determined) positive, square-free integers. 

\bcor
\label{COR:realquad}
If $K_*(C^*(R\rtimes R_+^\times))\cong K_*(C^*(R'\rtimes R_+'^\times))$, then $K=K'$. In particular, $K_*(C^*(R\rtimes R_+^\times))\cong K_*(C^*(R'\rtimes R_+'^\times))$ if and only if $C^*(R\rtimes R_+^\times)\cong C^*(R'\rtimes R_+'^\times)$ if and only if $K = K'$.
\ecor
\setlength{\parindent}{0cm} \setlength{\parskip}{0cm}

\bproof
We have $\#(R/(1-\epsilon)R)=|\det(1-\epsilon)|=\textup{trace}(\epsilon)-2$, so \Cref{thm:KCrealquad} allows us to retrieve $\textup{trace}(\epsilon)$ from the K-theory of our semigroup C*-algebra. Now let $D$ and $D'$ be the discriminants of $R_K$ and $R_{K'}$, respectively, and let $(t,u)$ and $(t',u')$ be the minimal positive solutions to the Pell equations $x^2-Dy^2= 4$ and $x^2-D'y^2= 4$, respectively, so that $\epsilon=\frac{t+u\sqrt{D}}{2}$ and $\epsilon'=\frac{t'+u'\sqrt{D'}}{2}$.
The equality $\textup{trace}(\epsilon)=\textup{trace}(\epsilon')$ means that $t=t'$, which forces $Du^2=D'u'^2$. Hence, $K=\Qz(\sqrt{D})=\Qz(\sqrt{Du^2})=\Qz(\sqrt{D'u^{'2}})=\Qz(\sqrt{D'})=K'$.
\eproof
Actually, the last equivalence in Corollary~\ref{COR:realquad} generalizes to Galois extensions of $\Qz$ (see \Cref{ss:RecC*}).
\setlength{\parindent}{0cm} \setlength{\parskip}{0.5cm}

\subsection{K-theory for boundary quotients}
\label{ss:Kbq}
\setlength{\parindent}{0cm} \setlength{\parskip}{0cm}

\subsubsection{Duality Theorems}

Let $\cS = \supp(\mfm_0)$, and set $\bm{R}_{\cS} = \prod_{v \in \cS^c, \, v \nmid \infty} R_v$, $\Az_{\cS} = \prod'_{v \in \cS^c, \, v \nmid \infty} K_v$. Here the restricted product is taken with respect to the subrings $R_v \subseteq K_v$, i.e., almost all coordinates of elements in $\Az_{\cS}$ belong to $R_v$.
\setlength{\parindent}{0cm} \setlength{\parskip}{0.5cm}

Our goal is to prove the following generalization of \cite[Theorem~4.1]{CL}:
\btheo
\label{THM:Duality}
$C_0(\Az_{\infty}) \rtimes Q$ and $C(\bm{R}_{\cS}) \rtimes R$ are $M$-equivariantly Morita equivalent.
\etheo
\setlength{\parindent}{0cm} \setlength{\parskip}{0cm}

We refer to \cite[\S~3.4]{LiLu_v1} for the notion of Morita equivalence equivariant under semigroup actions.
\setlength{\parindent}{0cm} \setlength{\parskip}{0.5cm}

Before we explain the proof, let us first record the following consequences:
\bcor
The Morita equivalence from Theorem~\ref{THM:Duality} induces an $M$-equivariant Morita equivalence $C_0(\Az_{\infty}) \rtimes Q \rtimes \mu \sim_M C(\Az_{\cS}) \rtimes Q \rtimes \mu$, which in turn induces ($M$-equivariant) isomorphisms 
\begin{equation}
\label{e:Kmu_inf-fin}
K_*(C_0(\Az_{\infty}) \rtimes Q \rtimes \spkl{M'}) \cong K_*(C(\bm{R}_{\cS}) \rtimes R \rtimes^e M') \ {\rm for} \ {\rm every} \ {\rm subsemigroup} \ M' \subseteq M \ {\rm with} \ \mu \subseteq M'.
\end{equation}
\ecor

\bcor
We have $C_0(\Az_{\infty}) \rtimes Q \rtimes G \sim_M C(\Az_{\cS}) \rtimes Q \rtimes G$.
\ecor

Let us now prove Theorem~\ref{THM:Duality}. The idea is to apply \cite[Proposition~3.30]{LiLu_v1} following the comment immediately after that proposition. The key observation is as follows: Let $\chi = \prod_v \chi_v$ be a character on $\Az_K$ inducing identifications $\widehat{\Az_K} \cong \Az_K$, $\widehat{\Az_{\infty}} \cong \Az_{\infty}$ and $\widehat{K} \cong \Az_K / K$ as constructed in Tate's thesis (see \cite{Tate} and \cite[Chapter~XIV]{Lang}). Write $\chi_{\cS} = \prod_{v \in \cS^c} \chi_v$, so that $\chi_{\cS}$ is a character on $\Az_{\infty} \times \Az_{\cS}$. 
\blemma
\label{LEM:Dual}
There exists $a \in K\reg$ such that the character $\chi_{\cS} \cdot a: \: \Az_{\infty} \times \Az_{\cS} \to \Tz, \, \bm{x} \ma \chi_{\cS}(a \bm{x})$ and the corresponding pairing $(\Az_{\infty} \times \Az_{\cS}) \times (\Az_{\infty} \times \Az_{\cS}) \to \Tz, \, (\bm{x}, \bm{y}) \ma (\chi_{\cS} \cdot a)(\bm{x} \bm{y})$ induce identifications $\widehat{\Az_{\infty}} \cong \Az_{\infty}$ and $\widehat{Q} \cong (\Az_{\infty} \times \Az_{\cS}) / Q$.
\elemma
\setlength{\parindent}{0cm} \setlength{\parskip}{0cm}

\bproof
$Q \subseteq \Az_{\infty} \times \Az_{\cS}$ is cocompact: Choose a compact subset $C_{\infty} \subseteq \Az_{\infty}$ such that $\Az_{\infty} = C_{\infty} + R$. Then let $C = C_{\infty} \times \bm{R}_{\cS}$. Given $\bm{x} \in \Az_{\infty} \times \Az_{\cS}$, by Strong Approximation we can find $t \in Q$ such that $\bm{x} - t \in \Az_{\infty} \times \bm{R}_{\cS}$. Now find $s \in R$ such that $\bm{x} - t - s \in C$. This shows that $\Az_{\infty} \times \Az_{\cS} = C + Q$. 
\setlength{\parindent}{0cm} \setlength{\parskip}{0.5cm}

Now let $\check{Q} = \menge{\bm{x} \in \Az_{\infty} \times \Az_{\cS}}{\chi_{\cS}(\bm{x} t) = 1 \ {\rm for} \ {\rm all} \ t \in Q}$. Then $\check{Q} \cong ((\Az_{\infty} \times \Az_{\cS}) / Q) \, \widehat{ } \ $ is discrete (as a dual group of a compact group). Moreover, we have $Q \subseteq \check{Q}$ because for every $t \in Q$, we have $\chi_v(t) = 1$ whenever $v \in \cS$ as $t \in R_v$, so that $\chi_{\cS}(t) = \chi(t) = 1$. Now $\check{Q} / Q \subseteq (\Az_{\infty} \times \Az_{\cS}) / Q$ is compact. Hence, being discrete and compact, $\check{Q} / Q$ must be finite. Thus there exists a positive integer $N$ such that $N \cdot \check{Q} \subseteq Q$, so that $\check{Q} \subseteq N^{-1} \cdot Q \subseteq K$. Let $\fC_v = \menge{x \in K_v}{{\rm Tr}_v(xy) \in \Zz_p \ {\rm for} \ {\rm all} \ y \in R_v}$, where $v \mid p$. This is (the local version of) Dedekind's complementary module in the sense of \cite[Chapter III, \S~2, Definition (2.1)]{Neu}. Let us now prove that $\check{Q} = \menge{x \in K}{x \in \fC_v \ {\rm for} \ {\rm all} \ v \in \cS}$. Surely, if $x \in K$ satisfies $x \in \fC_v$ for all $v \in \cS$, then $\chi_v(x \cdot R_v) = 1$ for all $v \in \cS$, so that $\chi_v(x \cdot Q) = \chi_v(x \cdot R_v) = 1$ for all $v \in \cS$, and thus $\chi_{\cS}(x \cdot Q) = \chi(x \cdot Q) = 1$. This proves \an{$\supseteq$}. For \an{$\subseteq$}, take $x \in \check{Q}$. Then $\chi_{\cS}(x \cdot Q) = 1$ and $\chi(x \cdot Q) = 1$. It follows that $(\chi \cdot \chi_{\cS}^{-1})(x \cdot Q) = 1$. As $Q$ is dense in $\prod_{v \in \cS} R_v$, we conclude that $\chi_v(x \cdot R_v) = 1$ for all $v \in \cS$. But this implies that $x \in \fC_v$ for all $v \in \cS$ because $\chi_v(x \cdot R_v) = \chi_p({\rm Tr}_v(x \cdot R_v)) = 1$ if and only if ${\rm Tr}_v(x \cdot R_v) \subseteq \Zz_p$ (where $v \mid p$). This proves \an{$\subseteq$}.

Now choose $a \in K\reg$ such that $a \cdot R_v = \fC_v$ for all $v \in \cS$. Here we are using Strong Approximation, and that $\cS$ is always finite. Then we have $\chi_{\cS}(a \bm{x} t) = 1$ for all $t \in Q$ if and only if $a \bm{x} \in \check{Q}$ if and only if $\bm{x} \in K$ and $a \bm{x} \in \fC_v = a \cdot R_v$ for all $v \in \cS$ if and only if $\bm{x} \in K$ and $\bm{x} \in R_v$ for all $v \in \cS$ if and only if $\bm{x} \in Q$. Thus under the pairing induced by $\chi_{\cS} \cdot a$, $Q$ is dual to itself, and hence we obtain the desired identification $\widehat{Q} \cong (\Az_{\infty} \times \Az_{\cS}) / Q$. The identification $\widehat{\Az_{\infty}} \cong \Az_{\infty}$ is already given by our choice of $\chi$.
\eproof

\bproof[Proof of Theorem~\ref{THM:Duality}]
Using Lemma~\ref{LEM:Dual}, the proof now proceeds as in \cite[\S~4]{CL}. First, observe that
$$
  C_0(\Az_{\infty}) \rtimes Q \cong C^*(\Az_{\infty}) \rtimes Q \cong C^*(Q) \rtimes \Az_{\infty} \cong C(\widehat{Q}) \rtimes \Az_{\infty} \cong C((\Az_{\infty} \times \Az_{\cS})/Q) \rtimes \Az_{\infty},
$$
where we used Lemma~\ref{LEM:Dual} in the last step. The (inverse of the) canonical multiplicative $M$-action on $C_0(\Az_{\infty}) \rtimes Q$ corresponds to the $M$-action $a. [t \ma f_t] = [t \ma \abs{N(a)}_{\infty}^{-1} f_{a^{-1}t} (a^{-1} \sqcup)]$ on $C((\Az_{\infty} \times \Az_{\cS})/Q) \rtimes \Az_{\infty}$. The second step is to prove that $C_c(\cG_N)$, equipped with analogous inner products and $M$-action as in \cite{LiLu_v1}, is an $M$-equivariant pre-imprimitivity bimodule in the sense of \cite{LiLu_v1}, which induces an $M$-equivariant $C((\Az_{\infty} \times \Az_{\cS})/Q) \rtimes \Az_{\infty}$-$C(\bm{R}_{\cS}) \rtimes R$-imprimitivity bimodule. Here $\cG$ is the transformation groupoid $((\Az_{\infty} \times \Az_{\cS})/Q) \rtimes \Az_{\infty}$ and $N$ is the image of $\gekl{0} \times \bm{R}_{\cS} \subseteq \Az_{\infty} \times \Az_{\cS}$ under the canonical projection $\Az_{\infty} \times \Az_{\cS} \onto (\Az_{\infty} \times \Az_{\cS})/Q$. These two steps together yield the desired $M$-equivariant Morita equivalence $C_0(\Az_{\infty}) \rtimes Q \sim_M C(\bm{R}_{\cS}) \rtimes R$.
\eproof
\setlength{\parindent}{0cm} \setlength{\parskip}{0.5cm}

\subsubsection{Rational K-theory computations for boundary quotients}

Using the duality theorem, we now compute the K-theory of the boundary quotient $\partial C^*_{\lambda}(R \rtimes M) \cong C(\bm{R}_{\cS}) \rtimes R \rtimes^e M$ rationally, under the assumption that 
\begin{equation}
\label{e:1-xi,m}
  \gcd \Big( \prod_{1 \neq \xi \in \mu} (1 - \xi), \ \mfm_0 \Big)=(1). 
\end{equation} 
Note that the boundary quotient corresponds to the maximal primitive ideal (see \S~\ref{ss:conmon} for the primitive ideal space computation), and its crossed product description can be derived as in \cite{Li14, EL} (see \cite[\S~8]{Bru}).

Here is the final result of our rational K-theory computation:
\btheo
\label{THM:QKtheory}
Assume that condition~\eqref{e:1-xi,m} holds. Then we have
\begin{eqnarray*}
  && \Qz \otimes K_*(\partial C^*_{\lambda}(R \rtimes M)) 
  \cong \Qz \otimes K_*(C(\bm{R}_{\cS}) \rtimes R \rtimes^e M)
  \cong \Qz \otimes K_*(C_0(\Az_{\cS}) \rtimes Q \rtimes G)\\
  &\cong& \Qz \otimes K_*(C_0(\Az_{\infty}) \rtimes Q \rtimes G)
  \cong \Qz \otimes K_*(C_0(\Az_{\infty}) \rtimes G), \ and
\end{eqnarray*}
$$
  \Qz \otimes K_*(C_0(\Az_{\infty}) \rtimes G)
  \cong
  \bfa
  \Qz \otimes \bigwedge^*(G) & \falls m=1,\\
  \Qz \otimes K_0(C^*(\mu)) \otimes \bigwedge^*(G/\mu) & \falls n \ {\rm even}, \, m>1,
  \efa
$$
if for every $c \in M$, the number $\# \menge{v_{\Rz}}{v_{\Rz}(c)<0}$ is even, whereas
$$
  \Qz \otimes K_*(C_0(\Az_{\infty}) \rtimes G)
  \cong
  \bfa
  \gekl{0} & \falls m=1,\\
  \Qz \otimes \bigwedge^*(G/\mu) & \falls m=2,
  \efa
$$
if there exists $c \in M$ with $\# \menge{v_{\Rz}}{v_{\Rz}(c)<0}$ odd.
\etheo
\setlength{\parindent}{0cm} \setlength{\parskip}{0cm}

Note that in all the cases where we get a non-zero answer, we always get $\bigoplus_{\Nz} \Qz$ as abstract vector spaces over $\Qz$. Moreover, the exterior algebras are equipped with their canonical gradings in all cases unless $n$ is odd, $m=1$, and for every $c \in M$, the number $\# \menge{v_{\Rz}}{v_{\Rz}(c)<0}$ is even. In that case, our computations yield an exterior algebra with reversed grading.
\bproof
The strategy is the one explained in \cite[Remark~3.16]{CL}, which is also used in \cite{LiLu}. First, using an inductive limit decomposition, compute (rational) K-theory for $C(\bm{R}_{\cS}) \rtimes R \rtimes \mu$. Then, for a convenient choice of $c \in M$, compute (rationally) the multiplicative action of $c$ on $C(\bm{R}_{\cS}) \rtimes R \rtimes \mu$ in K-theory. Next, use this computation together with Theorem~\ref{THM:Duality} (or rather \eqref{e:Kmu_inf-fin}) to show that the canonical inclusion $C_0(\Az_{\infty}) \rtimes G \into C_0(\Az_{\infty}) \rtimes Q \rtimes G$ induces a rational isomorphism. Finally, compute (rational) K-theory for $C_0(\Az_{\infty}) \rtimes G$ using a homotopy argument.
\setlength{\parindent}{0cm} \setlength{\parskip}{0.5cm}

In the following, we explain how to carry out each of these steps and summarize the final results. First, as in \cite[\S~4.2]{LiLu}, we find a decomposition $K_0(C^*(R \rtimes \mu)) = K_{\rm inf} \oplus K_{\rm fin}^c \oplus K_{\rm fin}^{\mu}$ such that $K_0(C(\bm{R}_{\cS}) \rtimes R \rtimes \mu) \cong \ilim_M \gekl{K_0(C^*(R \rtimes \mu), \eta_c}$ and for general $c \in M$ with $\prod_{1 \neq \xi \in \mu} (1 - \xi) \mid c$, $\id_{\Qz} \otimes \eta_c$ is of the form
$$
  \rukl{
  \begin{array}{c|c|c}
  A_c & * & 0 \\
  \hline
  0 & 0 & 0 \\
  \hline
  0 & 0 & \id
  \end{array}
  },
$$
whereas for $c \in \Zz_{>1}$, $c \in M$ with $\prod_{1 \neq \xi \in \mu} (1 - \xi) \mid c$, $\id_{\Qz} \otimes \eta_c$ is of the form
$$
  \rukl{
  \begin{array}{ccc|c|c}
  c^n & & 0 & & \\
   & \ddots & & * & 0 \\
  0 & & c^? & & \\
  \hline
   & 0 & & 0 & 0 \\
  \hline  
   & 0 & & 0 & \id
  \end{array}
  }
$$
with respect to the above decomposition of $K_0$ and suitable bases. Here $A_c: \: K_{\rm inf} \to K_{\rm inf}$ is an isomorphism. Moreover, the exponents of $c$ on the diagonal in the upper left box are non-negative and decreasing, and all the exponents are positive if $n$ is odd, while there is exactly one exponent equal to $0$ if $n$ is even.

For $K_1$, we also have $K_1(C(\bm{R}_{\cS}) \rtimes R \rtimes \mu) \cong \ilim_M \gekl{K_1(C^*(R \rtimes \mu), \theta_c}$. It follows from Theorem~\ref{LL} that $K_1$ vanishes if $m$ is even. If $m$ is odd, we know that $\id_{\Qz} \otimes \theta_c$ is an isomorphism for all $c \in M$, and for $c \in M$ with $c \in \Zz_{>1}$, $\id_{\Qz} \otimes \theta_c$ is of the form 
$$
  \rukl{
  \begin{array}{ccc}
  c^? & & 0 \\
   & \ddots & \\
  0 & & c^?
  \end{array}
  }
$$
with respect to a suitable basis, where again the exponents of $c$ on the diagonal are non-negative and decreasing, and there is exactly one exponent equal to $0$ if $n$ is odd, while all the exponents are positive if $n$ is even. 

It follows that $\Qz \otimes K_0(C(\bm{R}_{\cS}) \rtimes R \rtimes \mu) \cong \Qz^{r_0} \oplus \Qz^{m-1}$, and with respect to this decomposition, the multiplicative action of $c \in M$ with $c \in \Zz_{>1}$ and $\prod_{1 \neq \xi \in \mu} (1 - \xi) \mid c$ is given on $K_0$ by 
$$
  \beta^{\rm fin}_c
  =
  \rukl{
  \begin{array}{ccc|c}
  c^n & & 0 & \\
   & \ddots & & 0 \\
  0 & & c^? & \\
  \hline
   & 0 & & \id
  \end{array}
  }.
$$
For $K_1$, we obtain $\Qz \otimes K_1(C(\bm{R}_{\cS}) \rtimes R \rtimes \mu) \cong \Qz^{r_1}$, and the multiplicative action of $c \in M$ with $c \in \Zz_{>1}$ is given on $K_1$ by
$$
  \gamma_c^{\rm fin}
  = 
  \rukl{
  \begin{array}{ccc}
  c^? & & 0 \\
   & \ddots & \\
  0 & & c^?
  \end{array}
  },
$$
the same matrix from above describing $\id_{\Qz} \otimes \theta_c$. Hence the subgroup of $\Qz \otimes K_*(C(\bm{R}_{\cS}) \rtimes R \rtimes \mu)$ which is left invariant under $(\beta_c^{\rm fin},\gamma_c^{\rm fin})$ is precisely given by a one-dimensional subspace of $\Qz \otimes K_1$ (i.e., $\Qz \subseteq \Qz \otimes K_1$) if $n$ is odd and $m=1$, a one-dimensional subspace of $\Qz \otimes K_0$ (i.e., $\Qz \subseteq \Qz \otimes K_0$) if $n$ is even and $m=2$, and an $m$-dimensional subspace of $\Qz \otimes K_0$ (i.e., $\Qz^m \subseteq \Qz \otimes K_0$) in all other cases. 

Because of the equivariant K-theory identification in \eqref{e:Kmu_inf-fin}, the same description is valid for the multiplicative action corresponding to our element $c$ on $\Qz \otimes K_*(C_0(\Az_{\infty}) \rtimes Q \rtimes \mu)$. Comparing with the computation
$$
  \Qz \otimes K_*(C_0(\Az_{\infty}) \rtimes \mu) 
  \cong
  \bfa
  (\gekl{0},\Qz) & \falls n \ {\rm odd}, \, m = 1,\\
  (\Qz,\gekl{0}) & \falls n \ {\rm odd}, \, m = 2,\\
  (\Qz^m,\gekl{0}) & \sonst,
  \efa
$$
we see that the canonical inclusion $C_0(\Az_{\infty}) \rtimes \mu \into C_0(\Az_{\infty}) \rtimes Q \rtimes \mu$ induces a rational isomorphism onto the $\spkl{c}$-invariant part of $\Qz \otimes K_*(C_0(\Az_{\infty}) \rtimes Q \rtimes \mu)$. Extending $c \in M$ to a $\Zz$-basis $\gekl{c, c_1, c_2, \dotsc}$ of $G / \mu$, an iterative application of the Pimsner-Voiculescu exact sequence yields that the canonical inclusion $C_0(\Az_{\infty}) \rtimes \mu \rtimes \spkl{c, c_1, \dotsc, c_i} \into C_0(\Az_{\infty}) \rtimes Q \rtimes \mu \rtimes \spkl{c, c_1, \dotsc, c_i}$ induces a rational isomorphism in K-theory.  

All in all, we obtain 
\begin{eqnarray*}
  && \Qz \otimes K_*(\partial C^*_{\lambda}(R \rtimes M)) 
  \cong \Qz \otimes K_*(C(\bm{R}_{\cS}) \rtimes R \rtimes^e M)
  \cong \Qz \otimes K_*(C_0(\Az_{\cS}) \rtimes Q \rtimes G)\\
  &\cong& \Qz \otimes K_*(C_0(\Az_{\infty}) \rtimes Q \rtimes G)
  \cong \Qz \otimes K_*(C_0(\Az_{\infty}) \rtimes G).
\end{eqnarray*}

We now complete the proof by computing
$$
  \Qz \otimes K_*(C_0(\Az_{\infty}) \rtimes G)
  \cong
  \bfa
  \Qz \otimes \bigwedge^*(G) & \falls m=1,\\
  \Qz \otimes K_0(C^*(\mu)) \otimes \bigwedge^*(G/\mu) & \falls n \ {\rm even}, \, m>1,
  \efa
$$
if for every $c \in M$, the number $\# \menge{v_{\Rz}}{v_{\Rz}(c)<0}$ is even, and 
$$
  \Qz \otimes K_*(C_0(\Az_{\infty}) \rtimes G)
  \cong
  \bfa
  \gekl{0} & \falls m=1,\\
  \Qz \otimes \bigwedge^*(G/\mu) & \falls m=2,
  \efa
$$
if there exists $c \in M$ with $\# \menge{v_{\Rz}}{v_{\Rz}(c)<0}$ odd.
\eproof
\setlength{\parindent}{0cm} \setlength{\parskip}{0.5cm}

Let us now present an example showing why we only carry out rational computations.
\bex
Consider the case $K = \Qz$, $\mfm_0 = (2)$, $\mfm_{\infty}(\infty) = 1$, and $\Gamma = \gekl{+1} \times (\Zz / (2))^*$. Then $M = \menge{c \in \Zz_{>0}}{2 \nmid c}$, $Q = \menge{c \in \Qz_{>0}}{v_2(c) \geq 0}$ and $G = \menge{c \in \Qz_{>0}}{v_2(c) = 0}$. A similar computation as in \cite[\S~5]{Cun} and \cite[\S~5.1]{CL} yields
$$
  K_{\bullet}(C(\bm{R}_{\cS}) \rtimes R) \cong
  \bfa
  Q & \falls \bullet = 0,\\
  \Zz & \falls \bullet = 1.
  \efa
$$
Moreover, for $c=3$, the multiplicative action corresponding to $c$ is given by $3 \cdot \id_Q$ on $K_0$ and by $\id$ on $K_1$. Hence we obtain
$$
  K_{\bullet}(C(\bm{R}_{\cS}) \rtimes R \rtimes^e [c\rangle) \cong
  \bfa
  (\Zz / 2\Zz) \oplus \Zz & \falls \bullet = 0,\\
  \Zz & \falls \bullet = 1.
  \efa
$$
This example illustrates one of the reasons why we only carry out rational computations: Along the way of our computations, it could happen that torsion appears, which causes complications. Another reason is that the connecting maps in the inductive limit decomposition of $C(\bm{R}_{\cS}) \rtimes R \rtimes \mu$ are difficult to determine if we do not work over $\Qz$.
\eex

\bremark
However, it is possible to carry out precise K-theory computations for boundary quotients in some cases. For instance, assume condition~\eqref{e:1-xi,m} holds and that we can find $a, c \in M \cap \Zz_{>1}$ such that $c-a = 1$. Then we can compute K-theory without tensoring with $\Qz$. The point is that if we can find such elements $a$ and $c$, then the K-groups will be free abelian, so that they are completely determined by our rational computations. \eremark

Here is an example which explains why we need assumption \eqref{e:1-xi,m}.
\bex
Consider the case $K = \Qz$, $\mfm_0 = (2)$, $\mfm_{\infty}(\infty) = 0$, and $\Gamma = (\Zz / (2))^*$. Then $M = \menge{c \in \Zz\reg}{2 \nmid c}$, $Q = \menge{c \in \Qz\reg}{v_2(c) \geq 0}$ and $G = \menge{c \in \Qz\reg}{v_2(c) = 0}$. A similar computation as in \cite[\S~7]{Cun} and \cite[\S~3.1]{CL} yields
$$
  K_{\bullet}(C(\bm{R}_{\cS}) \rtimes R \rtimes \mu) \cong
  \bfa
  Q \oplus \Zz \oplus \Zz & \falls \bullet = 0,\\
  \gekl{0} & \falls \bullet = 1.
  \efa
$$
Moreover, for $c=3$, the multiplicative action corresponding to $c$ is given by the matrix
$$
  \begin{pmatrix}
  3 & 1 & 1 \\
  0 & 1 & 0 \\
  0 & 0 & 1
  \end{pmatrix}.
$$
Plugging this into the Pimsner-Voiculescu exact sequence, we obtain
$$
  K_{\bullet}(C(\bm{R}_{\cS}) \rtimes R \rtimes^e [c\rangle) \cong \Zz \oplus \Zz
$$
for $\bullet = 0, 1$. At the same time, we have $\Az_{\infty} = \Rz$ and $K_0(C_0(\Rz) \rtimes \mu) \cong \Zz$ and $K_1(C_0(\Rz) \rtimes \mu) \cong \gekl{0}$. This shows that the canonical embedding $C_0(\Rz) \rtimes \mu \into C_0(\Rz) \rtimes Q \rtimes \mu$ does not induce a rational isomorphism onto the part fixed by the multiplicative $\spkl{c}$-action. This also shows that for any $c'$ such that $\gekl{c,c'}$ can be extended to a $\Zz$-basis of $G/\mu$, the canonical embedding $C_0(\Rz) \rtimes \mu \rtimes \spkl{c,c'} \into C_0(\Rz) \rtimes Q \rtimes \mu \rtimes \spkl{c,c'}$ does not induce a rational isomorphism in K-theory.

This example shows why we need condition~\eqref{e:1-xi,m}. Otherwise we cannot compute (rational) K-theory of the boundary quotient by computing (rational) K-theory of $C_0(\Az_{\infty}) \rtimes G$.
\eex

\bremark
\label{rem:lbq--rbq}
As remarked in \cite[\S~5.11]{CELY}, it is an intriguing phenomenon that for cancellative semigroups, the left and right semigroup C*-algebras often have isomorphic K-theory. Actually, we do not know of an example where this does not happen because this phenomenon appears in all cases where we can compute K-theory. However, in general left and right semigroup C*-algebras have very different structural properties as C*-algebras (see \cite[\S~5.11]{CELY}). In the following, we give examples where the left and right boundary quotients have different K-theories:

Let $K$ be a number field with $r > 0$ real embeddings, choose $\mfm_{\infty}$ and $\Lambda$ so that $\mu = \gekl{+1}$, further choose $\mfm_0$ such that there exists $c \in M$ with $\# \menge{v_{\Rz}}{v_{\Rz}(c)<0}$ odd. Note that condition~\eqref{e:1-xi,m} is vacuous because $\mu$ is trivial. Now Theorem~\ref{THM:QKtheory} implies that we have 
$$
  \Qz \otimes K_*(\partial C^*_{\lambda}(R \rtimes M)) \cong \gekl{0}.
$$
However, for the right boundary quotient, we get $\partial C^*_{\rho}(R \rtimes M) \cong C^*(Q \rtimes G)$, and hence a straightforward computation shows that
$$
  \Qz \otimes K_*(\partial C^*_{\rho}(R \rtimes M)) 
  \cong \Qz \otimes K_*(C^*(Q \rtimes G)) 
  \cong \Qz \otimes K_*(C^*(G)) 
  \cong \Qz \otimes (\bigwedge{\hspace*{-0.1cm}}^* \, G)
  \cong \bigoplus_{\Nz} \Qz.
$$
In particular,
$$
  K_*(\partial C^*_{\lambda}(R \rtimes M)) \ncong K_*(\partial C^*_{\rho}(R \rtimes M)).
$$
Here is a concrete example: Let $K = \Qz[\sqrt{2}]$, so that $R = \Zz[\sqrt{2}]$. Let $v_+$ and $v_-$ be the real embeddings of $K$ determined by $v_{\pm}(\sqrt{2}) = \pm \sqrt{2}$. Set $\mfm_{\infty}(v_+) = 0$ and $\mfm_{\infty}(v_-) = 1$. Let $\Gamma = \gekl{+1}$ be trivial. Then $\mu = \gekl{+1}$ is trivial. Let $\mfm_0 = R$ be trivial. Then $c = 1 - \sqrt{2}$ lies in $M$; it satisfies $\# \menge{v_{\Rz}}{v_{\Rz}(c)<0} = \# \gekl{v_+} = 1$. 
\eremark

\section{Reconstruction theorems}

\subsection{Reconstruction using semigroup C*-algebras}
\label{ss:RecC*}

As discussed in \S~\ref{ss:conmon}, the non-zero minimal primitive ideals of $C^*_{\lambda}(R \rtimes M)$ correspond bijectively to the primes in $\cP_K^\m$. For each $\p\in\cP_K^{\mfm}$, we shall denote by $I_\p$ the minimal primitive ideal corresponding to $\p$.

It is straightforward to carry over the proof of \cite[Proposition~4.8]{Li14} to our more general situation, so that we obtain the following result.
\bprop[generalized version of {\cite[Proposition~4.8]{Li14}}]\label{prop:torsionorder}
For every $\p\in\cP_K^\m$ such that the canonical map $M^*\to (R/\p)^*$ is injective on $\mu$, the $K_0$-class $[1_{C^*_{\lambda}(R \rtimes M)/I_\p}]$ of the unit in $C^*_{\lambda}(R \rtimes M)/I_\p$ is a torsion element of order $\frac{N(\p)^{f(\p)}-1}{m}$.
\eprop

Let $p_{\max}:=\{p\in\Zz_{>0} \text{ prime } : p\mid N(1-\zeta^i) \text{ for some } 1\leq i\leq m-1\}$. For a prime $p\in\Zz_{>0}$, let $g_K(p)=\{\p\in\cP_K :\p\cap\Zz=p\Zz\}$ be the splitting number of $p$ in $K$. We collect the following consequences of Proposition~\ref{prop:torsionorder}.
\bcor
\label{cor:torsionorder}
$ $
\setlength{\parindent}{0cm} \setlength{\parskip}{0cm}

\begin{enumerate}
\item[(i)] For every prime $p\in\Zz_{>0}$ with $\frac{p-1}{m}>p_{\max}^{n \cdot \# C}-1$, we have
\[
g_K(p)=\#\{ I\in\Prim_{\min}(C^*_{\lambda}(R \rtimes M)) : p\mid (m\cdot \ord([1_{C^*_{\lambda}(R \rtimes M)/I}])+1)\}
\]
where $\Prim_{\min}(C^*_{\lambda}(R \rtimes M))$ denotes the set of non-zero minimal primitive ideals of $C^*_{\lambda}(R \rtimes M)$.
\item[(ii)] We have the following formula for the degree of $K$:
\[
[K:\Qz]=\limsup_{T\to\infty}\#\{I\in \Prim_{\min}(C^*_{\lambda}(R \rtimes M)) : \ord(C^*_{\lambda}(R \rtimes M)/I)=T\}.
\]
\end{enumerate}
\ecor
\setlength{\parindent}{0cm} \setlength{\parskip}{0cm}

\bproof
(i) and (ii) are generalized versions of \cite[Lemma~4.9]{Li14} and \cite[Lemma~5.1]{Li14}, respectively, and the proofs in \cite{Li14} carry over. 
\eproof
\setlength{\parindent}{0cm} \setlength{\parskip}{0.5cm}

Let $K$, $R$, $\mfm$, $\mu$, and $\Gamma$, $M$ be as above. Let $\fX \defeq {\rm Prim}_{\min}$ be the set of minimal non-zero primitive ideals of $C^*_{\lambda}(R \rtimes M)$. Define a function $o$ on $ \fX$ by setting $o(I) \defeq {\rm ord}(1_{C^*_r(R \rtimes M) / I})$. Let $n = [K:\Qz]$. Let $\rk = \rk(\Qz \otimes K_0(C^*_{\lambda}(R \rtimes M)))$ and set $e \defeq n + \rk$. Let $\cP$ denote the set of prime numbers (i.e., $\cP = \cP_{\Qz}$). 

The following allows us to read off the number of roots of unity:
\btheo
\label{thm:rootsofunity}
$m = \# \mu$ is the unique positive integer $m$ for which there exist functions $\fX \to \cP, \, I \ma p_I$ and $\fX \to \gekl{1, \dotsc, e}, \, I \ma i_I$ such that $I \ma p_I$ is finite-to-one and $\cP \setminus \menge{p_I}{I \in \fX}$ is finite, and we have
$$
  m \cdot o(I) + 1 = p_I^{i_I} \text{ for almost all } I \in \fX.
$$
\etheo
\bproof
It follows immediately from Proposition~\ref{prop:torsionorder} and Corollary~\ref{COR:rk=hrk} that $m$ has the desired properties. Note that Corollary~\ref{COR:rk=hrk} ensures that the codomain of $i$ is finite.

All we have to do is to prove uniqueness. Suppose $m'$ is another positive integer with $m \neq m'$ for which we can also find functions $\fX \to \cP, \, I \ma q_I$ and $\fX \to \gekl{1, \dotsc, e}, \, I \ma j_I$ with the above properties. Then we have for almost all $I \in \fX$:
\begin{equation}
\label{e:mp=mq+?}
  \frac{p_I^{i_I} - 1}{m} = \frac{q_I^{j_I} - 1}{m'} \ \Rarr \ m' p_I^{i_I} = m q_I^{j_I} + (m'-m).
\end{equation}
For every $j \in \gekl{1, \dotsc, e}$, consider $f_j(q) = mq^j + (m'-m)$ as a polynomial in $q$. For every non-constant polynomial $f$ with integer coefficients, there exist infinitely many primes $p$ for which there exists $z \in \Zz$ such that $f(z) \equiv 0 \  ({\rm mod} \ p)$ (see for instance \cite[Chapter~III, Theorem~45]{Nag}). Thus there exist pairwise distinct primes $p_j$, $1 \leq j \leq e$, and integers $z_j$, $1 \leq j \leq e$, such that, for all $1 \leq j \leq e$, $p_j \nmid m'$, $p_j \nmid m$, $p_j \nmid (m'-m)$, and $f_j(z_j) \equiv \ 0 \ ({\rm mod} \ p_j)$. In particular, we must have $p_j \nmid z_j$ for all $1 \leq j \leq e$. Here we used that $m \neq m'$, so that $m'-m \neq 0$.

Now set $N \defeq p_1 \dotsm p_e$ and find $z \in \Zz$ with $z \equiv z_j \ ({\rm mod} \ p_j)$ for all $1 \leq j \leq e$. Such $z$ exists by the Chinese Remainder Theorem. As $\gcd(z_j,p_j) = 1$ for all $1 \leq j \leq e$, we must have $\gcd(z,N) = 1$. Hence Dirichlet's Prime Number Theorem (see for instance \cite[Chapter~VII, (5.14)]{Neu}) implies that $\cQ \defeq \menge{q \in \cP}{q \equiv z \ ({\rm mod} \ N)}$ is infinite. As $\menge{q_I}{I \in \fX}$ contains almost all primes, $\menge{I \in \fX}{q_I \in \cQ}$ must be infinite. By \eqref{e:mp=mq+?}, we have
$$
  m' p_I^{i_I} = m q_I^{j_I} + (m'-m) \equiv 0 \ ({\rm mod} \ p_{j_I})
$$
for almost all $I \in \fX$ with $q_I \in \cQ$. As $p_{j_I} \nmid m'$, this implies $p_{j_I} \mid p_I^{i_I}$ and hence $p_I = p_{j_I}$ for almost all $I \in \fX$ with $q_I \in \cQ$. But $\menge{I \in \fX}{q_I \in \cQ}$ is infinite whereas $\menge{p_{j_I}}{I \in \fX_{\cQ}} \subseteq \menge{p_j}{1 \leq j \leq e}$ is finite, while $I \ma p_I$ is finite-to-one. This is a contradiction.
\eproof

\bcor\label{cor:rootsofunity}
Let $L$ be a another number field with ring of algebraic integers $S$,  and data $\mfn$, $\Lambda$ as above giving rise to the congruence monoid $N$. Let $\nu$ be the set of roots of unity in $N$. If $C^*_{\lambda}(R \rtimes M) \cong C^*_{\lambda}(S \rtimes N)$, then we must have $\# \mu = \# \nu$.
\ecor
\setlength{\parindent}{0cm} \setlength{\parskip}{0cm}
This answers the open question from \cite{Li14} whether it is possible to read off the number of roots of unity from our semigroup C*-algebras in the affirmative.
\setlength{\parindent}{0cm} \setlength{\parskip}{0.5cm}

Our next result shows that we can recover both the Dedekind zeta function of $K$ and the Kronecker set of $\Kz$ from K-theoretic invariants of $C_\lambda^*(R\rtimes M)$ (see \cite{Per} for the definition of arithmetic equivalence and its formulation in terms of Dedekind zeta functions and \cite{Jeh} for the definition of Kronecker equivalence). Under some additional assumptions on the number-theoretic input for our construction, this allows us to recover information about the congruence monoid $M$, the class field $\Kz$, and the class group $C$.

\btheo\label{thm:reconstruction}
Suppose that $K$ and $L$ are number fields with rings of algebraic integers $R$ and $S$. Let $\m$ and $\n$ be moduli for $K$ and $L$, and let $\Gamma$ and $\Lambda$ be subgroups of $(R/\m)^*$ and $(S/\n)^*$, respectively. Suppose that there is an isomorphism $C_\lambda^*(R\rtimes R_{\m,\Gamma})\cong C_\lambda^*(S\rtimes S_{\n,\Lambda})$. Then 
\setlength{\parindent}{0cm} \setlength{\parskip}{0cm}
\begin{enumerate}
\item[(i)] $K$ and $L$ are arithmetically equivalent, and $K(\m)^{\bar{\Gamma}}$ and $L(\n)^{\bar{\Lambda}}$ are Kronecker equivalent;
\item[(ii)] if the class fields $\Kmg$ and $L(\n)^{\bar{\Lambda}}$ of $K$ and $L$ are both Galois over $\Qz$, then $\#\clmg=\#\Cl_{\n}^{\bar{\Lambda}}$;
\item[(iii)] if $K$ or $L$ is Galois, then $K=L$; in particular, $K$ is Galois if and only if $L$ is Galois;
\item[(iv)] if $K$ or $L$ is Galois and both the class fields $\Kmg$ and $L(\n)^{\bar{\Lambda}}$ are Galois over $\Qz$, then
\begin{enumerate}
\item[(a)] $K=L$;
\item[(b)] $K(\m)^{\bar{\Gamma}}=L(\n)^{\bar{\Lambda}}$ (in any algebraically closed field containing both $K(\m)^{\bar{\Gamma}}$ and $L(\n)^{\bar{\Lambda}}$);
\item[(c)] $R^*\cdot(R_\n\cap R_{\m,\Gamma})=S^*\cdot (S_\m\cap S_{\n,\Lambda})$;
\item[(d)] $\clmg\cong\Cl_{\n}^{\bar{\Lambda}}$ (as abelian groups).
\end{enumerate}
\end{enumerate}
\etheo
\setlength{\parindent}{0cm} \setlength{\parskip}{0cm}

Here (and throughout our paper) we consider arithmetic equivalence and Kronecker equivalence over $\Qz$.
\bproof
(i): By Corollary~\ref{cor:rootsofunity}, $\rmg$ and $S_{\n,\Lambda}$ have the same number of roots of unity. Hence, combining \Cref{prop:torsionorder} and Corollary~\ref{cor:torsionorder}(i), we see that $g_K(p)=g_{L}(p)$ for all but finitely many primes. Therefore, $K$ and $L$ are arithmetically equivalent by \cite[Main~Theorem]{PerStu}. 
\setlength{\parindent}{0cm} \setlength{\parskip}{0.5cm}

Now let $\Kz=K(\m)^{\bar{\Gamma}}$ and $C=\clmg$. To show our claim about Kronecker equivalence, we need to show that we can recover, up to finitely many exceptions, the Kronecker set 
\[
D(\Kz\vert\Qz):=\{p\in\cP_\Qz : \exists \fP\in\cP_{\Kz}\text{ such that } \fP\mid p\text{ and } f_{\Kz/\Qz}(\fP\vert p)=1\}
\]
from K-theoretic invariants of $C_\lambda^*(R\rtimes R_{\m,\Gamma})$.
By Corollary~\ref{cor:rootsofunity}, we can read off the number of roots of unity in $\rmg$, so it follows from Proposition~\ref{prop:torsionorder} that we can recover the set $\{N(\p)^{f(\p)} : \p\in \cP_K^\m\}$ up to finitely many exceptions. Thus, we will be done if we show that $D(\Kz\vert\Qz)$ and $\cP_\Qz\cap \{N(\p)^{f(\p)} : \p\in \cP_K^\m\}$ differ by only finitely many elements. Our proof of this fact is inspired by \cite[Proof~of~Main~Theorem~(II)]{Glas}.
Let $\mathcal{R}:= \cP_\Qz \setminus \{p\in\cP_\Qz : p\notin\p \text{ for all }\p\in\supp(\m_0)\}$. Then $\mathcal{R}$ is finite, and every rational prime $p\in\mathcal{R}$ is unramified in $\Kz$ by Lemma~\ref{lem:splitsiff}. We will show that $D(\Kz\vert\Qz)\setminus\mathcal{R}=\cP_\Qz\cap \{N(\p)^{f(\p)} : \p\in \cP_K^\m\}\setminus\mathcal{R}$.

``$\subseteq$'': Suppose $p\in D(\Kz\vert\Qz)\setminus \mathcal{R}$. Since $N(\p)=p^{f_{K/\Qz}(\p\vert p)}$, we need to show that $f(\p)=1=f_{K/\Qz}(\p\vert p)$.
Choose any prime $\mathfrak{P}\in\cP_{\Kz}$ such that $\mathfrak{P} \mid p$ and $f_{\Kz/\Qz}(\mathfrak{P}\vert p)=1$, and let $\p:=\mathfrak{P}\cap R$. Then $\p\in\cP_K^\m$, and $\p$ lies over $p$. Since $f_{\Kz/\Qz}(\mathfrak{P}\vert p)=f_{\Kz/K}(\mathfrak{P}\vert \p)f_{K/\Qz}(\p\vert p)$, we see that $f_{\Kz/K}(\mathfrak{P}\vert \p)=f_{K/\Qz}(\p\vert p)=1$.
Since $\p\in\cP_K^\m$, $\p$ is unramified in $\Kz$ by Lemma~\ref{lem:splitsiff}. Thus, $f(\p)=f_{\Kz/\Qz}(\mathfrak{P}\vert \p)=1$ (see the proof of Lemma~\ref{lem:splitsiff}).

``$\supseteq$'': Now let $p\in\cP_\Qz\cap \{N(\p)^{f(\p)} : \p\in \cP_K^\m\}\setminus \mathcal{R}$. Then there exists $\p\in\cP_K^\m$ lying over $p$ with $f_{K/\Qz}(\p\vert p)=f(\p)=1$. Since $\p\in\cP_K^\m$, $\p$ is also unramified in $\Kz$, so $f(\p)=f_{\Kz/K}(\p)$. Now we have, for any prime $\fP$ of $\Kz$ lying over $\p$, that $f_{\Kz/\Qz}(\fP\vert p)=f_{\Kz/K}(\fP\vert\p)f_{K/\Qz}(\p\vert p)=f_{\Kz/K}(\p)f_{K/\Qz}(\p\vert p)=1$.

(ii): By part (i), we can recover the set $D(\Kz\vert\Qz)$, at least up to finitely many exceptions, from $C_\lambda^*(R\rtimes R_{\m,\Gamma})$. An application of the Chebotarev density theorem (\cite[Chapter~VII,~(13.4)~Theorem]{Neu}) shows that $D(\Kz\vert\Qz)$ has a Dirichlet density $\delta(D(\Kz\vert\Qz))$ (see, for instance, \cite[Chapter VII,~Definition~13.1]{Neu}) for the definition). Since $\Kz$ is Galois over $\Qz$, this Dirichlet density is given by $\delta(D(\Kz\vert\Qz))=\frac{1}{[K:\Qz]\cdot\#C}$ (see \cite[\S~3, equation~3.1]{Jeh}).
By Corollary~\ref{cor:torsionorder}(ii), we can read off $[K:\Qz]$ from K-theoretic invariants associated with $C_\lambda^*(R\rtimes R_{\m,\Gamma})$, so we can extract $\#C$.

(iii): By part (i), $K$ and $L$ are arithmetically equivalent, so $K$ and $L$ have the same Galois closure and degree by \cite[Theorem~1]{Per}. Hence, if $K$ or $L$ is Galois, so that it equals its Galois closure, then $K=L$. 

(iv): By part (iii), $K=L$, and by part (ii), the fields $K(\m)^{\bar{\Gamma}}$ and $K(\n)^{\bar{\Lambda}}$ are Kronecker equivalent. Since $K(\m)^{\bar{\Gamma}}$ and $K(\n)^{\bar{\Lambda}}$ are Galois over $\Qz$, a rational prime $p$ lies in the Kronecker set $D(K(\m)^{\bar{\Gamma}}\vert\Qz)$ if and only if $p$ splits completely in $K(\m)^{\bar{\Gamma}}$, and similarly for $K(\n)^{\bar{\Lambda}}$. This implies that a rational prime $p$ splits completely in $K(\m)^{\bar{\Gamma}}$ if and only if it splits completely $K(\n)^{\bar{\Lambda}}$. Therefore, $K(\m)^{\bar{\Gamma}}=K(\n)^{\bar{\Lambda}}$ in $\overline{\Qz}$ by \cite[Chapter~V,~Theorem~3.25]{MilCFT}. Now $R^*\cdot (R_\n\cap R_{\m,\Gamma})=R^*\cdot (R_\m\cap R_{\n,\Lambda})$ by Proposition~\ref{prop:monoidfromclassfield}, and we have $\clmg\cong\Gal(K(\m)^{\bar{\Gamma}}/K)=\Gal(K(\n)^{\bar{\Lambda}}/K)\cong \Cl_{\n}^{\bar{\Lambda}}$. 
\eproof
\setlength{\parindent}{0cm} \setlength{\parskip}{0.5cm}

\bremark
A necessary and sufficient condition for $\Kz$ to be Galois over $\Qz$ is given in Corollary~\ref{cor:kzgalois}. 
\eremark

In light of Theorem~\ref{thm:reconstruction}(iv), a natural question would be whether the existence of an isomorphism $C_\lambda^*(R\rtimes R_{\m,\Gamma})\cong C_\lambda^*(R\rtimes R_{\n,\Lambda})$ implies that $\supp(\m_0)=\supp(\n_0)$. The following partial answer is an immediate consequence of Proposition~\ref{prop:torsionorder}.

\blemma\label{lem:recoversupport}
Assume that $\mu_{\m,\Gamma}$ reduces injectively modulo $\p$ for all $\p\notin\supp(\m_0)$ and that $\mu_{\n,\Lambda}$ reduces injectively modulo $\p$ for every $\p\notin\supp(\n_0)$. Also assume that for every rational prime $p$, $\m_0$ is either divisible by every prime of $K$ over $\p$ or no primes of $K$ over $\p$, and similarly for $\n_0$.
If there is an isomorphism $C_\lambda^*(R\rtimes R_{\m,\Gamma})\cong C_\lambda^*(R\rtimes R_{\n,\Lambda})$, then $\supp(\m_0)=\supp(\n_0)$.
\elemma

\bremark
The hypothesis of Lemma~\ref{lem:recoversupport} is satisfied, for example, if $K$ is Galois, and $\sigma(\m_0)=\m_0$, $\sigma(\n_0)=\n_0$ for all $\sigma\in\Gal(K/\Qz)$, which is related to $\Kz$ being Galois over $\Qz$ (cf. Corollary~\ref{cor:kzgalois}).
\eremark

For the case $K=\Qz$, the above results can be strengthened. Indeed, we have the following result. 

\btheo
\label{thm:Q}
Let $\m$ and $\n$ be moduli for $\Qz$, and let $\Gamma$ and $\Lambda$ be subgroups of $(\Zz/\m)^*$ and $(\Zz/\n)^*$, respectively. If there is an isomorphism $C_\lambda^*(\Zz\rtimes\Zz_{\m,\Gamma})\cong C_\lambda^*(\Zz\rtimes\Zz_{\n,\Lambda})$, then 
\setlength{\parindent}{0cm} \setlength{\parskip}{0cm}
\begin{enumerate}
\item[(i)] $\Qz(\m)^{\bar{\Gamma}}=\Qz(\n)^{\bar{\Lambda}}$ (equality in $\overline{\Qz}$) and $\Cl_\m^{\bar{\Gamma}}\cong\Cl_\n^{\bar{\Lambda}}$;
\item[(ii)] if $-1\notin\Zz_{\m,\Gamma}$, then $-1\notin\Zz_{\n,\Lambda}$, $\langle\pm 1\rangle\cdot \Zz_{\m,\Gamma}=\langle\pm 1\rangle\cdot\Zz_{\n,\Lambda}$, and $\Zz_{\m,\Gamma}\subseteq \Zz_{>0}$ if and only if $\Zz_{\n,\Lambda}\subseteq \Zz_{>0}$.
\item[(iii)] if $-1\in\Zz_{\m,\Gamma}$, then $-1\in\Zz_{\n,\Lambda}$ and 
\begin{itemize}
\item $\Zz_{\m,\Gamma}=\Zz_{\n,\Lambda}$ provided that either $2\mid\m_0$ and $2\mid\n_0$, or $2\nmid\m_0$ and $2\nmid \n_0$;
\item $\{a\in \Zz_{\m,\Gamma} : \gcd(a,2q)=1\}=\{a\in \Zz_{\n,\Lambda}: \gcd(a,2q)=1\}$ for some odd prime $q$ if $2$ divides exactly one of $\m_0$ and $\n_0$.
\end{itemize}
\end{enumerate}
\etheo
\bproof
(i): This follows immediately from Theorem~\ref{thm:reconstruction}(iv) since $\Qz(\m)^{\bar{\Gamma}}$ and $\Qz(\n)^{\bar{\Lambda}}$ are Galois over $\Qz$.
\setlength{\parindent}{0cm} \setlength{\parskip}{0.5cm}

(ii): If $-1\notin\Zz_{\m,\Gamma}$, then it follows from Corollary~\ref{cor:rootsofunity} that $-1\notin\Zz_{\n,\Lambda}$ (another way to see this is to note that $K_1(C_\lambda^*(\Zz\rtimes\Zz_{\m,\Gamma}))$ vanishes if and only if $-1\in\Zz_{\m,\Gamma}$, and similarly for $\Zz_{\n,\Lambda}$, see \S~\ref{sss:K=Q}). 

Since $-1\notin \Zz_{\m,\Gamma}\cup\Zz_{\n,\Lambda}$, Lemma~\ref{lem:recoversupport} implies that $\supp(\m_0)=\supp(\n_0)$. Hence, Theorem~\ref{thm:reconstruction}(iv) implies that $\langle\pm1\rangle\Zz_{\m,\Gamma}=\langle\pm1\rangle\Zz_{\n,\Lambda}$.

Now Theorem~\ref{THM:QKtheory} implies that $\Qz\otimes K_*(\partial C_\lambda^*(\Zz\rtimes\Zz_{\m,\Gamma}))\neq\{0\}$ if and only if $\Zz_{\m,\Gamma}\subseteq \Zz_{>0}$ (and similarly for $\Zz_{\n,\Lambda}$). Thus, $\Zz_{\m,\Gamma}\subseteq \Zz_{>0}$ if and only if $\Zz_{\n,\Lambda}\subseteq \Zz_{>0}$.

(iii): If $-1\in \Zz_{\m,\Gamma}$, then it follows from Corollary~\ref{cor:rootsofunity} that $-1\in\Zz_{\n,\Lambda}$. We have three cases to consider.

Suppose $2\mid\m_0$ and $2\mid\n_0$. Since $\langle\pm1\rangle$ reduces injectively modulo $p$ for every odd prime $p$, Lemma~\ref{lem:recoversupport} implies that $\supp(\m_0)=\supp(\n_0)$, so that Theorem~\ref{thm:reconstruction}(iv) implies that $\langle\pm1\rangle\Zz_{\m,\Gamma}=\langle\pm1\rangle\Zz_{\n,\Lambda}$. Since $-1\in \Zz_{\m,\Gamma}\cap \Zz_{\n,\Lambda}$, it follows that $\Zz_{\m,\Gamma}=\Zz_{\n,\Lambda}$.

From the primitive ideal space computation discussed in \S~\ref{ss:conmon}, we see that the isomorphism $C_\lambda^*(\Zz\rtimes\Zz_{\m,\Gamma})\cong C_\lambda^*(\Zz\rtimes\Zz_{\n,\Lambda})$ induces a bijection $\varphi:\cP_\Qz^\m\cong \cP_\Qz^\n$ such that $C_\lambda^*(\Zz\rtimes\Zz_{\m,\Gamma})/I_p\cong C_\lambda^*(\Zz\rtimes\Zz_{\n,\Lambda})/J_{\varphi(p)}$ for all $p\in \cP_\Qz^\m$ where $I_p$ is the minimal primitive ideal of $C_\lambda^*(\Zz\rtimes\Zz_{\m,\Gamma})$ corresponding to $p\in \cP_\Qz^\m$ and $J_{\varphi(p)}$ is the minimal primitive ideal of $C_\lambda^*(\Zz\rtimes\Zz_{\n,\Lambda})$ corresponding to $\varphi(p)\in \cP_\Qz^\n$.

Now suppose $2\nmid\m_0$ and $2\nmid\n_0$, so that $2\in \cP_\Qz^\m\cap \cP_\Qz^\n$. Since $I_2$ contains the projection 
\[
1-\sum_{k\in\Zz/(2^{f^\Gamma_2})}\lambda(k,2^{f^\Gamma_2})\lambda(k,2^{f^\Gamma_2})^*
\] 
(see \cite[Theorem~7.1]{Bru}), it follows that $\ord([1_{C_{\lambda}^*(\Zz\rtimes\Zz_{\m,\Gamma})/I_2}]_0)$ divides $2^{f^\Gamma_2}-1$. Arguing as in the proof of \cite[Proposition~4.8]{Li14} and using analogues of \cite[Lemmas~4.3~\&~4.5]{Li14} shows that $\frac{2^{f^\Gamma_2}-1}{\gcd(2,2^{f^\Gamma_2}-1)}=2^{f^\Gamma_2}-1$ divides $\ord([1_{C_{\lambda}^*(\Zz\rtimes\Zz_{\m,\Gamma})/I_2}]_0)$, and thus $\ord([1_{C_{\lambda}^*(\Zz\rtimes\Zz_{\m,\Gamma})/I_2}]_0)=2^{f^\Gamma_2}-1$. Similarly, we have $\ord([1_{C_{\lambda}^*(\Zz\rtimes\Zz_{\n,\Lambda})/J_2}])=2^{f^\Lambda_2}-1$. There are now two sub-cases to consider.

If $\varphi(2)=q$ for some odd prime $q\in\cP_{\Qz}^\n$, and $\varphi^{-1}(2)=q'$ for some odd prime $q'\in \cP_\Qz^\m$, then, from the above discussion and Proposition~\ref{prop:torsionorder}, $2^{f^\Gamma_2}-1=\frac{q^{f^\Lambda_q}-1}{2}$ and $2^{f^\Lambda_2}-1=\frac{{q'}^{f^\Gamma_{q'}}-1}{2}$. 
Using (i), we have $f^\Gamma_2=\ord((2,\Qz(\m)^{\bar{\Gamma}}))=\ord((2,\Qz(\n)^{\bar{\Lambda}}))=f^\Lambda_2$, so we must have $q=q'$. Now $\varphi$ restricts to a bijection $\cP_\Qz^\m\setminus\{2,q\}\cong\cP_\Qz^\n\setminus\{2,q\}$ which must be the identity by Proposition~\ref{prop:torsionorder}; hence, $\supp(\m_0)\cup\{2,q\}=\supp(\n_0)\cup\{2,q\}$. Since $2,q\in\cP_\Qz^\m\cap\cP_\Qz^\n$, this implies $\supp(\m_0)=\supp(\n_0)$. As in the first case above, Theorem~\ref{thm:reconstruction}(iv) and the fact that $-1\in \Zz_{\m,\Gamma}\cap \Zz_{\n,\Lambda}$ imply that $\Zz_{\m,\Gamma}=\Zz_{\n,\Lambda}$.

If $\varphi(2)=2$, then Proposition~\ref{prop:torsionorder} implies that $\varphi$ must be the identity map which forces $\supp(\m_0)=\supp(\n_0)$, so that Theorem~\ref{thm:reconstruction}(iv) implies that $\langle\pm1\rangle\Zz_{\m,\Gamma}=\langle\pm1\rangle\Zz_{\n,\Lambda}$. Since $-1\in \Zz_{\m,\Gamma}\cap \Zz_{\n,\Lambda}$, it follows that $\Zz_{\m,\Gamma}=\Zz_{\n,\Lambda}$.

For the last case, suppose without loss of generality that $2\nmid \m_0$ and $2\mid \n_0$. Let $q\in\cP_{\Qz}^\n$ be the odd prime such that $\varphi(2)=q$. Then, by Proposition~\ref{prop:torsionorder}, $\varphi$ restricts to the identity map from $\cP_\Qz^\m\setminus\{2\}$ onto $\cP_\Qz^\n\setminus\{q\}$, which implies that $\supp(\m_0)\cup\{2\}=\supp(\n_0)\cup\{q\}$. Now Theorem~\ref{thm:reconstruction}(iv) implies that $\langle\pm 1\rangle\cdot(\{a\in \Zz_{\m,\Gamma} : \gcd(a,2q)=1\})= \langle\pm 1\rangle\cdot(\{a\in \Zz_{\n,\Lambda}: \gcd(a,2q)=1\})$. Since $-1$ is in $\{a\in \Zz_{\m,\Gamma} : \gcd(a,2q)=1\}$ and $\{a\in \Zz_{\n,\Lambda}: \gcd(a,2q)=1\}$, we are done. \eproof
\setlength{\parindent}{0cm} \setlength{\parskip}{0.5cm}

\subsection{Reconstruction using Cartan pairs}
\label{ss:RecCartan}

As discussed in \S~\ref{ss:conmon}, let $D_{\lambda}(R \rtimes M)$ denote the canonical Cartan subalgebra of $C^*_{\lambda}(R \rtimes M)$.
\btheo
\label{thm:Cartan}
Let $K$, $L$ be number fields with rings of algebraic integers $R$, $S$, and suppose that we are given data $\mfm$, $\Gamma$ and $\mfn$, $\Lambda$ as in \S~\ref{ss:conmon} for $K$ and $L$, respectively. Let $M$ and $N$ be the corresponding congruence monoids. If $(C^*_{\lambda}(R \rtimes M), D_{\lambda}(R \rtimes M)) \cong (C^*_{\lambda}(S \rtimes N), D_{\lambda}(S \rtimes N))$, then there is a bijection $\varphi: \: \cP^{\mfm}_K \cong \cP^{\mfn}_L$ such that $N(\mfp) = N(\varphi(\mfp))$, $f^{\Gamma}_{\mfp} = f^{\Lambda}_{\varphi(\mfp)}$, and we have ${\rm Cl}_{\mfm}^{\bar{\Gamma}} \cong {\rm Cl}_{\mfn}^{\bar{\Lambda}}$ (as abelian groups).
\etheo
\setlength{\parindent}{0cm} \setlength{\parskip}{0cm}

\bproof
In the following, let us fix a number field $K$ with ring of algebraic integers $R$ and given data $\mfm$, $\Gamma$, and let us explain how to recover the prime ideals of $K$ that are relatively prime to $\m_0$, together with the functions $N(-)$, $f(-)$, as well as the group $C \defeq {\rm Cl}_{\mfm}^{\bar{\Gamma}}$, from the Cartan pair $(C^*_{\lambda}(R \rtimes M), D_{\lambda}(R \rtimes M))$. We proceed as in \cite{Li16_2}.
\setlength{\parindent}{0cm} \setlength{\parskip}{0.5cm}

To simplify notation, since the number field and all the relevant data are fixed we drop sub- and superscripts following our notational conventions (see \S~\ref{s:pi}). Write $D \defeq D_{\lambda}(R \rtimes M)$. For every subset $\cF$ of $\cP = \cP^{\mfm}_K$, let $I_{\cF}$ be the primitive ideal of $C^*_{\lambda}(R \rtimes M)$ corresponding to $\cF$, set $D_{\cF} \defeq I_{\cF} \cap D$, and let $\iota_{\cF}: \ D_{\cF} \into D$, $i: D \into C^*_{\lambda}(R \rtimes M)$ be the canonical embeddings. We denote the induced homomorphisms in $K_0$ by $(\iota_{\cF})_*$ and $i_*$. Let $\Delta \defeq i_*(K_0(D))$, $\Delta_{\cF} \defeq i_*((\iota_{\cF})_*(K_0(D_{\cF})))$ and write $\pi_{\cF}$ for the canonical projection $\Delta \onto \Delta / \Delta_{\cF}$. For $\cF = \emptyset$, we set $D_{\cF} \defeq (0)$ and $\Delta_{\cF} = \gekl{0}$. Given a collection $\cF$ of prime ideals in $\cP$, let $C_{\cF}$ be the subgroup of $C$ given by $\spkl{\menge{[\mfp]}{\mfp \in \cF}} \subseteq C$, where $C_{\emptyset}$ is the trivial subgroup.

It follows from the K-theory formula in Theorem~\ref{thm:ktheoryformula} that $\Delta \cong \bigoplus_C \Zz$. Let $M_{\mfp}$ be the composite $\Delta \overset{\mfp_*}{\lori} \Delta \overset{N(\mfp) \, \id}{\lori} \Delta$, $[e_{\mfa}] \ma N(\mfp) [e_{\mfp \mfa}]$, where $e_{\mfa}$ denotes the canonical projection in $D \subseteq C^*_{\lambda}(R \rtimes M)$ corresponding to $\mfa$. Using the observation (proven as in \cite{Li16_2})
\begin{equation}
\label{e:D_F}
  \Delta_{\cF} = \sum_{\mfp \in \cF} (\id - M_{\mfp})(\Delta) \ \text{for every subset} \ \cF \subseteq \cP,
\end{equation}
applied to singletons $\cF = \gekl{\mfp}$, we obtain as an application of \cite[Lemma~2.3]{Li16_2} that
\begin{equation}
\label{e:D/D_p}
  \Delta / \Delta_{\gekl{\mfp}} \cong \bigoplus_{C / \spkl{[\mfp]}} \Zz / (N(\mfp)^{\# \spkl{[\mfp]}} - 1) \Zz \ \text{for every} \ \mfp \in \cP.
\end{equation}
Moreover, \eqref{e:D_F} and \cite[Lemma~2.3]{Li16_2} imply the following technical result (the proof is as in \cite{Li16_2}):
\bprop[analogue of Proposition~2.1 in \cite{Li16_2}]
\label{PROP:D/D_F}
Let $\cF$ be a finite collection of prime ideals $\mfp$ in $\cP$ with $\mfp \nmid 2$.
\setlength{\parindent}{0.5cm} \setlength{\parskip}{0cm}

There exists $d_{\cF} \in \Nz_0$ with $\Delta / \Delta_{\cF} \cong \bigoplus_{C / C_{\cF}} \Zz / d_{\cF} \Zz$. For ${\cF} = \emptyset$, we have $d_{\emptyset} = 0$, and for ${\cF} \neq \emptyset$, $d_{\cF}$ is positive and even.

Moreover, there are ideals $\mfa_{{\cF},i}$ in $\cP$ with $C = \dotcup_i \, C_{\cF} [\mfa_{{\cF},i}]$ and such that $\gekl{\pi_{\cF} [e_{\mfa_{{\cF},i}}]}_i$ forms a $(\Zz / d_{\cF} \Zz)$-basis for $\Delta / \Delta_{\cF}$. Given an ideal $\mfa$ in $\cP$ with $[\mfa] \in C_{\cF} [\mfa_{{\cF},i}]$, there exists an odd number $l_{\cF}(a) \in \Nz$ with $\pi_{\cF} [e_{\mfa}] = l_{\cF}(\mfa) \pi_{\cF} [e_{\mfa_{{\cF},i}}]$.
\eprop
\setlength{\parindent}{0cm} \setlength{\parskip}{0cm}

Now we complete the proof of Theorem~\ref{thm:Cartan} as follows: We can read off $\cP$ (as a set) as the minimal non-zero primitive ideals of $C^*_{\lambda}(R \rtimes M)$. We can further read of $\# C$ as $\# C = \rk_{\Zz}(\Delta)$, and thus \eqref{e:D/D_p} allows us to reconstruct the functions $N(-)$ as well as $f(-)$ on $\cP$. Finally, Proposition~\ref{PROP:D/D_F} together with \cite[Lemma~3.2]{Li16_2} allows us to reconstruct the group $C$.
\eproof
\setlength{\parindent}{0cm} \setlength{\parskip}{0.5cm}


As an immediate consequence, we obtain that in the case of the rational number field, isomorphism of Cartan pairs yields a stronger conclusion compared to the one in Theorem~\ref{thm:Q}.
\bcor
\label{cor:CartanQ}
Let $\m$ and $\n$ be moduli for $\Qz$, and let $\Gamma$ and $\Lambda$ be subgroups of $(\Zz/\m)^*$ and $(\Zz/\n)^*$, respectively. If there is an isomorphism $(C_\lambda^*(\Zz\rtimes\Zz_{\m,\Gamma}),D_\lambda(\Zz\rtimes\Zz_{\m,\Gamma}))\cong (C_\lambda^*(\Zz\rtimes\Zz_{\n,\Lambda}),D_\lambda(\Zz\rtimes\Zz_{\n,\Lambda}))$, then $\langle\pm1\rangle\Zz_{\m,\Gamma}=\langle\pm1\rangle\Zz_{\n,\Lambda}$.
\ecor
\setlength{\parindent}{0cm} \setlength{\parskip}{0cm}
\bproof
By Theorem~\ref{thm:Cartan}, there exists a bijection $\varphi: \: \cP^{\mfm}_\Qz \cong \cP^{\mfn}_\Qz$ such that $N(p) = N(\varphi(p))$. Hence, $\varphi$ must be the identity map, so $\supp(\mfm_0)=\supp(\mfn_0)$. Now, as in the proof of Theorem~\ref{thm:Q}, Theorem~\ref{thm:reconstruction}(iv) implies that $\langle\pm1\rangle\Zz_{\m,\Gamma}=\langle\pm1\rangle\Zz_{\n,\Lambda}$.
\eproof

\setlength{\parindent}{0cm} \setlength{\parskip}{0.5cm}
Let us combine Theorem~\ref{thm:Cartan} with the following number-theoretic result.
\bprop
Assume that $[\Kmg:K]=[L(\n)^{\bar{\Lambda}}:L]$ and that there is a bijection $\varphi:\cP_K^\m\to \cP_L^{\n}$ such that $N(\p)=N(\varphi(\p))$ and $f_\p^\Gamma=f_{\varphi(\p)}^\Lambda$ for all $\p\in\cP_K^\m$. Then $K$ and $L$ are arithmetically equivalent, and $\Kmg$ and $L(\n)^{\bar{\Lambda}}$ are arithmetically equivalent.
\eprop
\setlength{\parindent}{0cm} \setlength{\parskip}{0cm}

\bproof
For all but finitely many rational primes $p$, we have
\[
g_{K/\Qz}(p)=\#\{\p\in \cP_K^\m : p\mid N(\p)\}=\#\{\varphi(\p)\in \cP_L^{\n} : p\mid N(\varphi(\p))\}=g_{L/\Qz}(p).
\]
Hence, $K$ and $L$ are arithmetically equivalent by \cite[Main~Theorem]{PerStu}.
\setlength{\parindent}{0cm} \setlength{\parskip}{0.5cm}

Since every prime $\p\in \cP_K^\m$ is unramifed in $\Kmg$ by (b) in \S~\ref{sec:cft}, we have, for each $\p\in \cP_K^\m$, 
\[
g_{\Kmg/K}(\p)f_{\Kmg/K}(\p)=[\Kmg:K]
\] 
where $g_{\Kmg/K}(\p)$ denotes the splitting number of $\p$ in $\Kmg$, and similarly 
\[
g_{L(\n)^{\bar{\Lambda}}/L}(\q)f_{L(\n)^{\bar{\Lambda}}/L}(\q)=[L(\n)^{\bar{\Lambda}}:L]
\]
for all $\q\in\cP_L^\n$. Since $N(\p)=N(\varphi(\p))$, we see that, for a rational prime $p$, we have $\p\mid p$ if and only if $\varphi(\p)\mid p$. Thus, for all but finitely many rational primes $p$, $\varphi$ restricts to a bijection from the set of primes of $K$ lying over $p$ onto the set of primes of $L$ lying over $p$. That is, for all but finitely many rational primes $p$, if $pR=\prod_{i=1}^k\p_i$ is the prime factorization of $pR$ in $R$ with the $\p_i$ distinct primes, then the prime factorization of $pS$ in $S$ is given by $pS=\prod_{i=1}^k\varphi(\p_i)$, and the $\varphi(\p_i)$ are distinct. For any such prime $p$, we have
\begin{align*}
g_{\Kmg/\Qz}(p)&=\sum_{i=1}^kg_{\Kmg/K}(\p_i) 
=\sum_{i=1}^k\frac{[\Kmg:K]}{f_{\Kmg/K}(\p)}
=\sum_{i=1}^k\frac{[L(\n)^{\bar{\Lambda}}:L]}{f_{L(\n)^{\bar{\Lambda}}/L}(\varphi(\p))}
=\sum_{i=1}^kg_{L(\n)^{\bar{\Lambda}}/L}(\varphi(\p_i))\\
&=g_{L(\n)^{\bar{\Lambda}}/\Qz}(p).
\end{align*}
where the middle equality used our assumption that $[\Kmg:K]=[L(\n)^{\bar{\Lambda}}:L]$. Hence, $\Kmg$ and $L(\n)^{\bar{\Lambda}}$ are split equivalent, so $\Kmg$ and $L(\n)^{\bar{\Lambda}}$ are arithmetically equivalent by \cite[Main~Theorem]{PerStu}.
\eproof
\setlength{\parindent}{0cm} \setlength{\parskip}{0.5cm}

\bcor
\label{cor:Cartan}
Let $K$, $L$ be number fields with rings of algebraic integers $R$, $S$, and suppose we are given data $\mfm$, $\Gamma$ and $\mfn$, $\Lambda$ for $K$ and $L$, respectively. Let $M$ and $N$ be the corresponding congruence monoids. If $(C^*_{\lambda}(R \rtimes M), D_{\lambda}(R \rtimes M)) \cong (C^*_{\lambda}(S \rtimes N), D_{\lambda}(S \rtimes N))$, then $K$ and $L$ are arithmetically equivalent, $\Kmg$ and $L(\n)^{\bar{\Lambda}}$ are arithmetically equivalent and we have ${\rm Cl}_{\mfm}^{\bar{\Gamma}} \cong {\rm Cl}_{\mfn}^{\bar{\Lambda}}$.
\ecor

\bremark
We can reformulate Corollary~\ref{cor:Cartan} in terms of continuous orbit equivalence of partial dynamical systems as in \cite[Theorem~1.2]{Li16_2}.
\eremark

\end{document}